\documentclass[preprint,12pt]{elsarticleplain}

\usepackage{caption}
\usepackage{subfigure}
\usepackage{amssymb}
\usepackage{amsthm}
\usepackage{amsmath}
\usepackage{mathtools}
\usepackage{algorithm}
\usepackage{algpseudocode}
\usepackage[hidelinks]{hyperref} 
\usepackage{diagbox}
\usepackage{multirow}
\usepackage{algorithmicx}
\usepackage{tikz}
\usepackage{booktabs}
\usepackage{adjustbox}

\def\grad{\nabla}

\def\d{~\mathrm{d}}

\def\Id{\mathrm{Id}}
\DeclareMathOperator{\diag}{diag}

\def\bx{\boldsymbol{x}}
\def\bv{\boldsymbol{v}}

\newtheorem{theorem}{Theorem}[section]

\numberwithin{equation}{section}

\usepackage{xcolor}

\usepackage[normalem]{ulem}

\DeclareMathOperator*{\argmax}{arg\,max}
\algrenewcommand{\algorithmicrequire}{\textbf{Inputs:}}

\begin{document}

\begin{frontmatter}



\title{Adaptive and hybrid reduced order models to mitigate Kolmogorov barrier in a multiscale kinetic transport equation}


\author[inst1]{Tianyu Jin} 
\author[inst1]{Zhichao Peng\textsuperscript{*}}
\author[inst1,inst2]{Yang Xiang\textsuperscript{**}}

\cortext[*]{Corresponding author, pengzhic@ust.hk (Zhichao Peng)}
\cortext[**]{Corresponding author, maxiang@ust.hk (Yang Xiang)}

\affiliation[inst1]{organization={Department of Mathematics},
            addressline={The Hong Kong University of Science and Technology}, 
            country={Hong Kong SAR}}
\affiliation[inst2]{organization={Algorithms of Machine Learning and Autonomous Driving Research Lab},
            addressline={HKUST Shenzhen-Hong Kong Collaborative Innovation Research Institute}, 
            city={Futian},
            state={Shenzhen},
            country={China}}

\begin{abstract}
In this work, we develop reduced order models (ROMs) to predict solutions to a multiscale kinetic transport equation with a diffusion limit under the parametric setting. When the underlying scattering effect is not sufficiently strong, the system governed by this equation exhibits transport-dominated behavior. Suffering from the Kolmogorov barrier for transport dominant problems, classical linear ROMs may become inefficient in this regime. To address this issue, we first develop a piecewise linear ROM by introducing a novel goal-oriented adaptive time partitioning strategy. To avoid local over-refinement or under-refinement, we propose an adaptive coarsening and refinement strategy that remains robust with various initial empirical partitions. Additionally, for problems where a local linear approximation is not sufficiently efficient, we further develop a hybrid ROM, which combines autoencoder-based nonlinear ROMs and piecewise linear ROMs. Compared to previous autoencoder-based ROMs, this hybridized method reduces the offline autoencoder's training cost by only applying it to time intervals that are adaptively identified as the most challenging. Numerical experiments  demonstrate that our proposed approaches successfully predict full-order solutions at unseen parameter values with both efficiency and accuracy.
To the best of our knowledge, this is the first attempt to address the Kolmogorov barrier for multiscale kinetic transport problems with the coexistence of both transport- and diffusion-dominant behaviors.  
\end{abstract}




\end{frontmatter}



\section{Introduction}\label{sec:intro}
The radiative transfer equation (RTE) is a kinetic equation modeling particle systems propagating through and interacting with a background medium. This equation, whose unknown can be seen as a particle distribution function posed in a high-dimensional phase space, finds a wide range of applications across diverse scientific and engineering disciplines, including medical imaging, nuclear engineering, astrophysics and remote sensing. However, the inherent high-dimensional nature of RTE poses significant challenges for high-fidelity numerical simulations. 
This computational burden becomes especially pronounced in multi-query scenarios where repeated numerical solutions for various model parameters are required. As a result, effective dimensionality reduction techniques for RTE are highly desired to make computational costs manageable in real-time or multi-query applications.

Reduced order modeling (ROM) \cite{benner2015survey} has emerged as a powerful data-driven technique for constructing low-rank approximations to high-dimensional parametric problems. ROM typically employs an offline-online decomposition: during the offline phase, low-dimensional structures are extracted by exploring data for a parametric problem, while the online phase leverages these low-rank structures to rapidly predict solutions for new parameters via projection or interpolation. Recently, in the context of RTE, ROMs have been actively and rapidly advanced through various methodologies, including proper orthogonal decomposition (POD) \cite{buchan2015pod, choi2021space, tano2021affine, behne2022minimally, behne2023parametric, coale2023reduced, coale2024reduced,halvic2023non,hardy2024proper}, dynamical mode decomposition (DMD) \cite{mcclarren2022data,smith2023variable}, and greedy algorithms \cite{tencer2016reduced,peng2022reduced,peng2024micro}. Beyond  predictions, ROMs have have also been leveraged to accelerate linear solvers and eigenvalue solvers for RTE \cite{mcclarren2019calculating,roberts2019acceleration,mcclarren2022data,smith2023variable,peng2024reduced,peng2024flexible}.

Despite the initial success of ROMs for RTE, limitations persist in previous studies. Most existing ROMs employ a linear approximation that uses low-dimensional linear subspaces to approximate the underlying solution manifold. However, these linear methods face fundamental limitations when applied to transport-dominant problems, as characterized by the Kolmogorov barrier \cite{slowdecay,rim2023performance,peherstorfer2022breaking}. This barrier manifests a slow decay of the approximation error with increasing dimension of the reduced space -- a fundamental challenge for RTE systems where particle behavior becomes transport-dominant in regimes with insufficiently strong absorption and scattering effects. Although nonlinear model order reduction (MOR) techniques have emerged as potential solutions to overcome the Kolmogorov barrier in various transport problems, their applications to RTE have not been thoroughly investigated. This work seeks to bridge this gap in the literature by developing nonlinear MOR approaches specifically tailored for multiscale RTE.

Before introducing our methods, we first briefly review two predominant strategies for overcoming the Kolmogorov barrier for transport problems.  The first strategy is to utilize piecewise linear approximations to the solution manifold by introducing adaptivity. Examples following this idea include time partitioning \cite{adpative_time_partition, san2015principal,goal-oriented,ji2024aaroc,li2025localized}, local spatial refinement \cite{carlberg2015adaptive}, and online adaptive update of reduced space \cite{peherstorfer2020model}. The second strategy is to utilize nonlinear approximations based on intrinsic rank-reveal nonlinear transformations. Such nonlinear transformations can be constructed using characteristic line information 
\cite{rowley2000reconstruction,mojgani2017lagrangian, reiss2018shifted,lu2020lagrangian}, optimal transport \cite{ehrlacher2020nonlinear,rim2023manifold}, registration \cite{taddei2020registration}, and deep learning \cite{lee2020model,kim2022fast,papapicco2022neural,peng2023learning,barnett2023neural}. Instead of an exhaustive literature review, we direct readers to the review paper \cite{peherstorfer2020model} for a more complete overview.

In this paper, we first develop a piecewise linear ROM for RTE based on a goal-oriented adaptive time partitioning strategy. Unlike existing approaches \cite{adpative_time_partition} that begin with a single interval and only contain refinement, our method introduces two key innovations: (1) a coarsening strategy enabling robust use of arbitrary empirical initial partitions, and (2) an equilibrium detection step to allow an ultra-low dimensional reduced order space when the system approaches its long-time steady-state.

Although piecewise linear ROMs based on time partitioning have demonstrated effectiveness in capturing local low-rank structures of the solution manifold \cite{san2015principal,peherstorfer2020model}, these methods face limitations when the solution manifold undergoes sharp changes within a short time interval, such as the propagation of a delta function. In such cases, achieving sufficient accuracy through time refinement often necessitates excessive partitioning, resulting in prohibitive memory cost. For such cases, autoencoder-based nonlinear ROMs \cite{lee2020model,choi2021space} become an attractive choice as they offer superior capability in capturing low-rank structures compared to piecewise linear ROMs. However, their training phase is more expensive than linear or piecewise linear ROMs. To leverage advantages of both methods, we further propose a hybrid approach that strategically applies autoencoder-based nonlinear ROMs locally in time intervals where linear reduction is not sufficiently efficient.

Building upon the above discussions, we aim to develop computationally efficient and reliable nonlinear ROMs for multiscale problems governed by RTE. The main contributions of this work are summarized as follows.
\begin{itemize}
    \item We further develop the adaptive time partitioning strategy in \cite{adpative_time_partition} which only involves refinement by introducing coarsening and equilibrium detection. Coarsening allows us to use various initial partitions without worrying about unnecessary local refinement, while equilibrium detection improves the efficiency for long-time simulations.
    \item When the solution manifold possesses sharp features in a short time, piecewise linear ROMs may lead to excessive local refinement or require a large number of bases in some time intervals.
    To address this issue, we introduce an adaptive strategy to identify such intervals and employ an autoencoder-based nonlinear ROM in them. To the best of our knowledge, this is the first work to adaptively combine local nonlinear ROMs with linear ROMs. 
    
    \item The proposed methods outperform the traditional linear model reduction method for multiscale kinetic transport equations with transport-dominant subregions in terms of both efficiency and accuracy. Additionally, to the best of the authors' knowledge, this is the first nonlinear model reduction attempt to mitigate the Kolmogorov barrier in kinetic transport equations that encompass both parabolic and hyperbolic regions.
\end{itemize}
The remainder of this paper is organized as follows. In Section \ref{sec:setting}, we present the formulation of the problem and provide essential background information that will be pertinent in subsequent sections. Section \ref{sec:method} details the proposed methods, with a particular focus on offline strategies, and includes thorough analysis of reconstruction error and complexity. In Section \ref{sec:num_experiments}, we evaluate the efficiency and effectiveness of proposed methods through numerical examples exhibiting multiscale properties. Finally, we draw our conclusions in Section \ref{sec:conclusion}. 

\section{Problem setting and preliminaries}\label{sec:setting}

Throughout this manuscript we focus on solving the following time-dependent linear kinetic transport equation with single energy group and isotropic scattering:
\begin{align}\label{eq:kinetic}
\partial_t f + \bv\cdot\grad_{\bx} f = \sigma_s(\rho-f) - \sigma_a f + G, 
\end{align}
equipped with initial condition $f(\bx,\bv,t=0)=f_0$ and appropriate boundary conditions. Here, $f(\textbf{x},\textbf{v},t)$ is the particle distribution at spatial location $\textbf{x}\in\Omega_{\textbf{x}}\subset\mathbb{R}^d\hspace{2pt}(d=1,2)$, time $t\in\mathbb{R}^+$ with velocity $\textbf{v}\in\Omega_{\textbf{v}}$. The velocity space $\Omega_v$ is $[-1,1]$ for 1D slab geometry, and the unit sphere $\mathbb{S}^2$ for higher dimensions.
$\sigma_s(\textbf{x})\geq0$, $\sigma_a(\textbf{x})\geq0$ and $\sigma_t = \sigma_s+\sigma_a$ are the isotropic scattering cross section, absorption cross section and total cross section, respectively. The macroscopic density is defined as the average of the distribution function over the velocity space, i.e. $\rho(\bx,t)=\int_{\bv\in\Omega_v} f(\bx,\bv,t)d\nu$, where $\nu$ is a measure of the velocity space satisfying $\int_{\Omega_{\textbf{v}}}\d\nu=1$. When $\Omega_v=[-1,1]$, $\int_{\bv\in\Omega_v}f\d\nu=\frac{1}{2}\int_{-1}^1 fdv$. When $\Omega_v=\mathbb{S}^2$, $\int_{\bv\in\Omega_v}f\d\nu=\frac{1}{4\pi}\int_{\bv\in\mathbb{S}^2} fd\bv$. $G(\textbf{x},t)$ is an isotropic source term.

To discretize equation \eqref{eq:kinetic}, we apply the discrete ordinates ($S_N$) method \cite{pomraning1973equations} in the velocity space. Specifically, we solve equation \eqref{eq:kinetic} at a set of quadrature points for the velocity space $\{(\bv_j,\omega_j)\}_{j=1}^{N_{\bv}}$, where $\omega_j$'s are corresponding quadrature weights. When $\Omega_v=[-1,1]$, Gauss-Legendre points are utilized, while Chebyshev-Legendre points are considered when $\Omega_v=\mathbb{S}^2$.
The $S_N$ system of \eqref{eq:kinetic} is therefore given by
\begin{equation}\label{eq:SN}
\partial_t f_j +\bv_j\cdot\grad_{\textbf{x}} f_j = \sigma_s(\rho-f_j) - \sigma_a f_j + G,
\quad \rho=\sum_{j=1}^{N_{\bv}}\omega_j f_j.
\end{equation}
Here, $f_j$ is the approximation of the distribution at velocity $\bv_j$.
To fully discretize the $S_N$ system \eqref{eq:SN}, we further apply the linear upwind discontinuous Galerkin (DG) method in space \cite{adams2001discontinuous} along with the second-order backward differential formula (BDF2) in time. It is noteworthy that the linear upwind DG method is proved to be asymptotic preserving \cite{adams2001discontinuous,guermond2010asymptotic}.

\subsection{Reduced order model and proper orthogonal decomposition}
Due to the high-dimensional nature of equation \eqref{eq:kinetic}, numerical simulations for it is usually computationally expensive, especially for parametric problems arising from multi-query applications such as design optimization and uncertainty quantification. Reduced order model (ROM) \cite{benner2015survey}, as we mentioned in Section \ref{sec:intro}, is a dimensionality reduction technique, which leverages low-rank structures of the solution manifold for parametric problems. 

ROMs are typically constructed following an offline-online decomposition framework. In the offline stage, a low-dimensional approximation to the solution manifold is constructed by extracting low-rank structure of the solution manifold from data, i.e. high-fidelity numerical solutions. Then in the online stage, fast prediction or reconstruction is performed through interpolation or projection, leveraging the low-dimensional  approximation to the solution manifold constructed offline. 

In the following, we briefly review the idea of proper orthogonal decomposition and outline an interpolation-based non-intrusive online stage.

\subsubsection{Proper orthogonal decomposition}\label{sec:pod}
Proper orthogonal decomposition (POD) constructs a low-dimensional reduced order space by computing singular value decomposition (SVD) of a snapshot matrix.

Let the vector representing the high-fidelity solution corresponding to time $t_i$ and parameter $\mu_j$ be $f_{h}(t_i;\mu_j)\in\mathbb{R}^{n_h}$.  The snapshot matrix is defined as
\begin{align*}
    S = [f_h(t_1;\mu_1)-\bar{f}_h, \dots,f_h(t_{n_t};\mu_1)-\bar{f}_h, \dots,f_h(t_{n_t};\mu_{n_p})-\bar{f}_h]\in\mathbb{R}^{n_h\times n_s}.
\end{align*}
Here, $n_h = n_{\textbf{x}}\times n_{\textbf{v}}$ is the dimension of the high-fidelity solution, $n_s = n_t\times n_p$ is the total number of snapshots, $\bar{f}_h\in\mathbb{R}^{n_h}$ is a time-independent offset, which is the mean of the set $\{f_h(t_i;\mu_j):i=1,\dots,n_t; j=1,\dots,n_p\}$. Each column of the snapshot matrix is the difference between a high-fidelity solution and the average of high-fidelity solutions.
whose columns are vectors representing full-order solutions for the training set $\{(t_i,\mu_j): 1\leq i\leq n_t, 1\leq j\leq n_p\}$.
 
To build reduced order space, we start by computing the singular value decomposition (SVD) of the snapshot matrix,
\begin{align*}
    S = U\Sigma V^\top ,
\end{align*}
where diagonal elements of $\Sigma = \diag\{\sigma_1,\sigma_2,\dots,\sigma_s\}\in\mathbb{R}^{n_h\times n_s}$ are the singular values with $\sigma_1\geq\sigma_2\geq\dots\geq\sigma_s\geq0$ and $s \leq \min\{n_h,n_s\}$. The columns of $U\in\mathbb{R}^{n_h\times n_h}$ and $V\in\mathbb{R}^{n_s\times n_s}$ are the left and right singular vectors, respectively. POD method constructs a reduced order space $U_r\in\mathbb{R}^{n_h\times r}$ which is the first $r$ columns of the matrix $U$. Its dimension is determined by finding the minimum $r$ satisfying 
\begin{align}
\frac{\sigma_{r+1}}{\sigma_1}\leq \varepsilon_{\text{POD}},\label{eq:2-norm-truncation}
\end{align}
or
\begin{align}
\frac{\sum_{i=1}^r \sigma_i^2}{\sum_{i=1}^s \sigma_i^2}\geq 1-\varepsilon^2_{POD}.\label{eq:f-norm-truncation}
\end{align}
By the Schmidt-Eckart-Young theorem, the best rank-$r$ linear approximation of the snapshot matrix $S$, is then given by $U_r\Sigma_rV_r^T$, where $\Sigma_r=\textrm{diag}\{\sigma_1,\dots,\sigma_r\}\in\mathbb{R}^{r\times r}$ and $V_r$ is the first $r$ columns of $V$, and the relative truncation error according to \eqref{eq:2-norm-truncation} or \eqref{eq:f-norm-truncation} satisfies 
$||U_r\Sigma_rV_r^\top -S||/||S||\leq \varepsilon_{\textrm{POD}}$, where $||\cdot||$ denotes $||\cdot||_2$ or $||\cdot||_F$ correspondingly. 

In the online stage of POD, a reduced order solution satisfying $f_{h,r}=\bar{f}_h+U_rc_r$, where $c_r\in\mathbb{R}^r$ is the coordinates in the reduced order space, will be constructed through projection or interpolation.

\subsubsection{A non-intrusive online stage}\label{sec:online}
For simplicity, we apply a simple non-intrusive online stage to predict solutions for new parameters, as our main focus is the development of the offline stage. 

Our approach involves two primary steps: first, interpolating low-dimensional reduced order representations for new parameters, and second, reconstructing high-dimensional full-order solutions from these low-dimensional representations.

After constructing the reduced order space, we can derive the reduced order representation of solution as $f(t_i;\mu_j)\approx \bar{f}_h+U_rc_r(t_i;\mu_j)$ with $c_r\in\mathbb{R}^{r\times n_s}$. In the online stage, when predicting solutions to a new parameter value, we interpolate these low-dimensional reduced order representations corresponding to our training parameters
$$
c_r(t,\tilde{\mu}) = \mathcal{I}(t,\tilde{\mu}) = \mathcal{I}(c_r(t;\mu_1),\dots,c_r(t;\mu_{n_p})).
$$
In our implementation, when the dimension of the parameter $d_p$ is equal to $1$, we employ piecewise polynomial interpolation. For higher dimensional parameters, the radial basis function (RBF) interpolation is adopted:
\begin{align*}
     \mathcal{I}(t,\tilde{\mu}) =\sum_{i=1}^{n_p}w_i\varphi(\|\tilde{\mu}-\mu_i\|),
\end{align*}
where $\varphi(\|\mu-\mu_i\|)$ is a kernel function. The weights $w_i$'s are determined by enforcing
\begin{align*}
    \sum_{i=1}^{n_p}w_i\varphi(\|\mu_j-\mu_i\|) = c_r(t,\mu_j), \quad j=1,\dots, n_p.
\end{align*}
 In our experiments, we select a quintic kernel function.

Finally, after obtaining the reduced order representation, the prediction for the solution corresponding to the new parameter $\tilde{\mu}$ can be reconstructed through $f(t,\tilde{\mu})=\bar{f}_h+U_rc_r(t,\tilde{\mu})$. 

Although the online stage we presented is based on a linear reduced order space, it can be extended to accommodate nonlinear ROMs by replacing the linear combination of reduced bases with a nonlinear mapping from the reduced order representation to the full order solution.
 
\subsection{Kolmogorov barrier for ROM construction}\label{sec:kolmogorov-barrier}
Classical linear ROM such as POD utilizes an $r$-dimensional linear space, $\mathcal{F}_r$, to approximate the solution manifold, namely $\mathcal{M}_h$, of the parametric problem. The efficacy of linear ROMs relies on the fast decay of the Kolmogorov $n$-width 
\begin{equation}\label{eq:kol_n_width}
d_n(\mathcal{M}_h)=\inf _{\mathcal{F}_n, \operatorname{dim}\left(\mathcal{F}_n\right)=n} \sup _{u \in \mathcal{M}_h} \inf _{u_n \in \mathcal{F}_n}\|u-u_n\|,
\end{equation}
which is the optimal error  of approximating the solution manifold with a $n$-dimensional linear space.
However, as discussed in \cite{ohlberger2015reduced,peherstorfer2022breaking}, the Kolmogorov $n$-width may decay slowly for transport-dominant problems, necessitating a high-dimensional reduced order space to obtain an accurate reduced order approximation.

Unfortunately, when the scattering and absorption effects are weak, the behavior of the particle system modeled by equation \eqref{eq:kinetic} becomes transport-dominant. Particularly, in the free-streaming limit with zero scattering and absorption effects (i.e. $\sigma_s=\sigma_a=0$), equation \eqref{eq:kinetic} will be reduced to a purely transport equation. As a result, without sufficiently strong scattering or absorption effects, classical linear ROMs may become ineffective for \eqref{eq:kinetic}.

Nonlinear ROMs, which employ nonlinear approximations to the solution manifold, have been developed to mitigate the Kolmogorov barrier \cite{peherstorfer2022breaking}. However, compared to linear ROMs, the construction of non-linear ROMs typically incurs higher computational costs. Notably, for diffusion-dominant problems, linear ROMs may achieve a favorable balance between online efficiency and offline construction costs. Moreover, under conditions of sufficiently strong scattering effects or over long time-scale, the behavior of the particle system modeled by \eqref{eq:kinetic} transitions to the diffusion-dominant regime, making linear ROMs a compelling choice for such cases.

In the following sections, we will develop a hybrid ROM designed specifically for transport-dominant problems or multiscale problems governed by equation \eqref{eq:kinetic} with transport-dominant sub-regions. Our approach integrates autoencoder-based nonlinear approximations \cite{lee2020model} with piecewise linear in time approximations. The goal is to leverage linear approximations whenever they are sufficiently efficient while selectively employing autoencoders to avoid over-refinement in the time domain. To the best of our knowledge, this represents the first attempt to mitigate the Kolmogorov barrier in multiscale kinetic problems and combine piecewise linear in time methods with autoencoder-based ROMs.

\section{Offline time partitioning to mitigate the Kolmogorov barrier}\label{sec:method}
The ultimate goal of this section is to mitigate the Kolmogorov barrier in the development of efficient ROMs for kinetic transport problems that exhibit transport-dominant phenomena. 

Although the solution manifold of a transport-dominant parametric problem may not be well approximated by a global linear ROM, efficient piecewise linear ROMs can be built via time partitioning techniques, as discussed in previous studies \cite{adpative_time_partition,san2015principal}. Uniform partitioning \cite{san2015principal}, though easy to implement, often poses difficulties in determining the optimal number of intervals necessary to achieve the desired online efficiency without sacrificing accuracy \textit{a priori}. 
To overcome this limitation, we introduce a goal-oriented adaptive partitioning strategy that incorporates both coarsening and refinement capabilities, which aims to create robust and efficient ROMs tailored specifically for kinetic transport problems.

Despite the initial success of piecewise linear ROMs for various transport-dominant problems, they may lose efficiency when confronted with challenging problems that exhibit sharp features in time, e.g. the propagation of a delta function. For such problems, maintaining accuracy often necessitates excessive time refinement, which can lead to significant memory costs associated with storing local basis functions across all intervals. To avoid over-refinement, we enhance our time partitioning strategy by integrating advanced nonlinear ROMs that utilize autoencoder neural networks  \cite{lee2020model}. Although autoencoders offer superior compression capabilities compared to traditional linear or piecewise linear ROMs, their implementation incurs higher offline computational costs due to the training of the neural network. Also, due to the highly non-convex nature of neural networks, the optimizer may get stuck in a local minima during the training process. To effectively balance the trade-off between compression performance and offline costs while avoiding over-refinement, we strategically combine time partitioning techniques with autoencoder-based ROMs, deploying autoencoders selectively in regions where they are most beneficial.

In what follows, we start from reviewing uniform time partitioning in Section \ref{sec:uniform} to demonstrate basic ideas, then introduce our novel goal-oriented adaptive partitioning in Section \ref{sec:ATPCPOD} followed by the autoencoder-based hybrid ROM in Section \ref{sec:hybridPOD}.

\subsection{Principal interval decomposition (uniform partitioning)}\label{sec:uniform}
We begin with recalling uniform time partitioning to present essential ideas. Uniform time partitioning has been applied to the principal interval decomposition (PID) in \cite{san2015principal}. Although uniform partition may not be optimal, its simplicity and effectiveness have been demonstrated for various challenging convection-dominated problems, as presented in \cite{san2015principal}.

We uniformly divide the entire time domain $[0,T]$ into $k\in \mathbb{N}^+$ time intervals $\{\mathcal{T}_i = ((i-1)\Delta t, i\Delta t]\}_{i=1}^{k}$, where $\Delta t = \frac{T}{k}$. The snapshot matrix $S$ is therefore correspondingly divided as $S=(S_1|S_2|\dots|S_k)$, where $S_i$ represents snapshots that lie in the time interval $\mathcal{T}_i$. Then following Section \ref{sec:pod}, within each time interval $\mathcal{T}_i$, POD with a tolerance of $\varepsilon_{\text{POD}}$ is applied to construct a reduced basis represented by the orthogonal matrix $B_i\in\mathbb{R}^{n_h\times n_{r_i}}$. In each time interval, we seek the reduced order solution in the $r_i$-dimensional linear space determined by $B_i$. 

As shown in Theorem \ref{thm:error_frob}, the overall reconstruction error of the entire snapshot matrix remains bounded by $\varepsilon_{\textrm{POD}}$ regardless of the value of $k$ in the Frobenius norm. 
\begin{theorem}\label{thm:error_frob} (\textbf{Reconstruction error.})
    Given the snapshot matrix $S\in\mathbb{R}^{n_h\times n_s}$ and partition its columns into $k$ parts, then we obtain $k$ submatrices $S_j\in\mathbb{R}^{n_h\times n_j}, j=1,2,\dots,k$ and $\sum_{j=1}^k n_j = n_s$, i.e. $S = (S_1|\dots|S_k)$. Suppose for any $j = 1,2,\dots, k$, the best $r_j$-rank approximation for $S_j$ is $S_{r_j}$. 
    
    \noindent
    If the relative approximation error of $S_{r_j}$ satisfies
    $$
    \frac{\|S_j-S_{r_j}\|_F}{\|S_j\|_F}\leq\varepsilon \text{ for all }j=1,2,\dots,k,
    $$
    then the relative reconstruction error for the whole snapshot matrix satisfies
    $$
    \frac{\|S-(S_{r_1}|\dots|S_{r_k})\|_F}{\|S\|_F}\leq\varepsilon.
    $$
\end{theorem}
The proof of this theorem is straightforward, and we leave it to \ref{sec:proof_error_frob}. 
We emphasize that Theorem \ref{thm:error_frob} can be applied to arbitrary partitions.

Next, we analyze the online and offline complexity of classical POD on the whole time interval and the POD with uniform time partitioning, respectively.

\textbf{Online computational savings.} In the online stage, the computational cost of interpolating reduced order representation is proportional to the dimension of the reduced order space, $r$. If a projection-based online stage is applied, such a cost will typically be proportional to $r^3$ or $r^2$ depending on the time discretization. After uniform partitioning, as long as the dimension of reduced order space for each time interval, $r_i$, is smaller than the dimension of reduced order space using the classical POD method, $r$, the computational savings in the online stage would be expected. 

\textbf{Offline complexity.} Classical POD constructs reduced order space through SVD of the snapshot matrix $S\in\mathbb{R}^{n_h\times n_s}$. The computational cost associated with SVD, using Golub-Kahan-Reisch algorithm \cite{svd_scipy}, is $O(n_hn_s^2)$. After uniform partitioning, the number of snapshots in each interval becomes $\frac{n_s}{k}$. The total cost associated with SVD on all intervals becomes $k\times O(n_h (\frac{n_s}{k})^2)=O(n_h\frac{n_s^2}{k})$. To interpolate reduced order solution online, we need to construct reduced order representation for training parameters offline, which leads to $O(rn_hn_s)$ cost in classical POD and $O(n_h\frac{n_s}{k}(\sum_{i=1}^k r_i))$ cost after uniform partitioning with $r_i$ being the dimension of the reduced order space for time interval $\mathcal{T}_i$. 

After discussing the time complexity, we now examine the associated memory costs in offline stage. The offline memory usage can be split into two primary components: (i) the cost of storing basis matrices and (ii) the cost of storing low-dimensional coordinates in the latent spaces. For standard POD, the memory cost for storing basis matrices is $\mathcal{O}(n_h\times r)$, while after uniform partitioning it necessitates $\mathcal{O}(n_h\sum_{i=1}^k r_i)$ in memory. Regarding the storage of coordinates, classical POD requires $\mathcal{O}(r\times n_s)$ memory, while a memory usage of $\mathcal{O}(\frac{n_s}{k}\sum_{i=1}^k r_i)$ is demanded after uniform partitioning. It is important to note that in typical scenarios involving high-dimensional problems, $n_h\gg n_s$, making (i) the dominant storage component.
Though the dimension of the reduced order space in each time interval $r_i$ will typically be smaller than the dimension for the classical POD, $r$, the total number of basis $\sum_{i=1}^k r_i\geq r$, which incurs higher memory costs than classical POD. 
This is because the classical SVD provides an optimal rank-$r$ reconstruction to the snapshot matrix. Even using the same truncation tolerance as the classical POD, the reconstruction given by uniform partitioning, $(S_{r_1}|\dots|S_{r_k})$, is sub-optimal. 

\begin{figure}[htbp]
	\centering
       \subfigure[]{
		\begin{minipage}{0.45\textwidth}
			\includegraphics[width=1\textwidth]{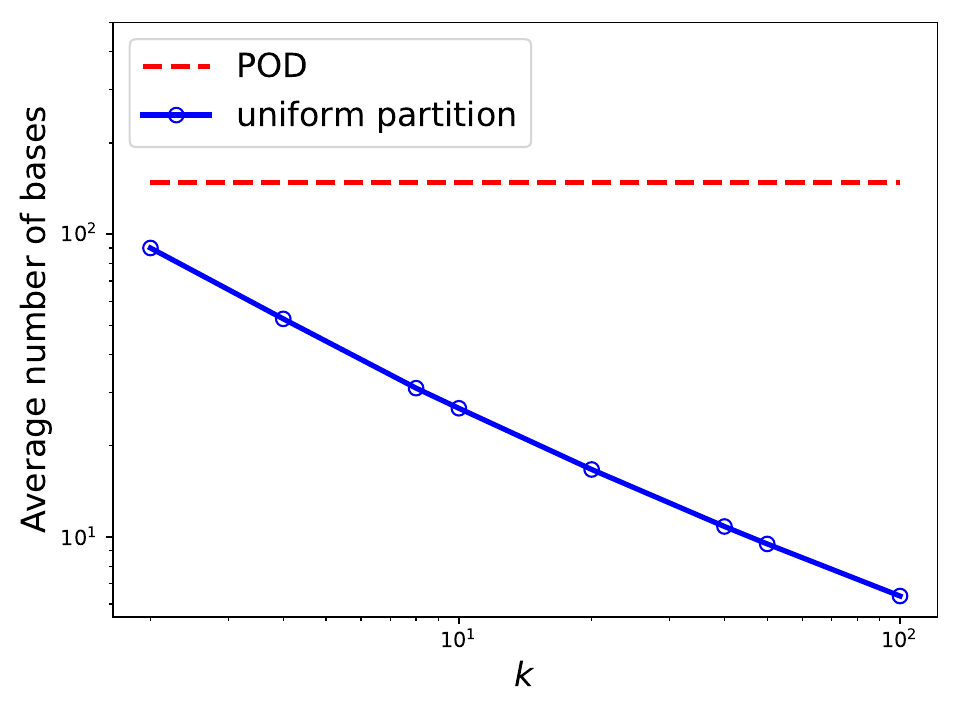}
		\end{minipage}
	}%
        \subfigure[]{
    		\begin{minipage}{0.45\textwidth}
   		 	\includegraphics[width=1\textwidth]{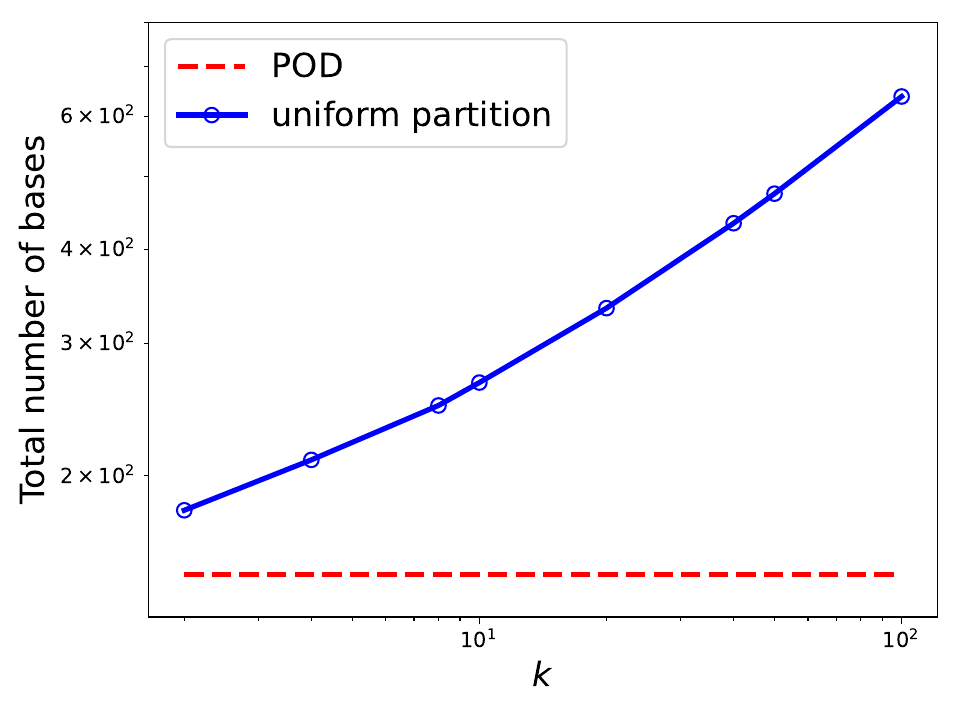}
    		\end{minipage}
    	}
	\caption{Complexity of classical POD and uniform time partition into $k$ intervals in Example 1 when varying $k$, $\varepsilon_{\text{POD}}=10^{-4}$. (a) Average number of bases regarding online time; (b) Total number of bases regarding memory cost.}
	\label{fig:complexity}
\end{figure}

The average and total number of bases of POD with uniform time partitioning (blue lines) are presented in Figure \ref{fig:complexity} compared to those of classical POD (red lines). The results perfectly corroborate our complexity analysis above. Specifically, Figure \ref{fig:complexity} (a) shows that a finer uniform time partitioning leads to enhanced computational efficiency during the online stage. Nevertheless, Figure \ref{fig:complexity} (b) reveals the price of the acceleration achieved: when the number of intervals, $k$, becomes excessively large, particularly in cases of over-refinement, the uniform time partitioning incurs substantial memory overhead.

In summary, paying the price of higher memory compared to classical POD, time partitioning achieves acceleration in the online stage by paying the memory price to save more reduced basis.

\subsection{Goal-oriented adaptive time partitioning}\label{sec:ATPCPOD}
The goal of introducing time partitioning is to achieve desired online efficiency without losing accuracy. When uniform partitioning is employed, determining the optimal number of partitions \textit{a priori} can be challenging and may lead to local under- or over-refinement. To address this issue, an adaptive refinement strategy is introduced in \cite{adpative_time_partition} under the POD-greedy framework. However, it always starts from one interval and lacks the coarsening ability.
Here, we introduce a goal-oriented adaptive partitioning strategy that incorporates both refinement and coarsening capabilities. This new strategy enables robust utilization of empirical initial partitions while mitigating risks of local over-refinement. Additionally, we introduce an equilibrium detection step to allow the use of an ultra-low-dimensional reduced order space when the underlying dynamical system approaches its equilibrium.

Our key idea is to refine time intervals with too many basis functions and merge intervals with too few basis functions to their neighbours. The detailed procedure is as follows.

Let $\hat{r}_{M}$ and $\hat{r}_{m}\leq \hat{r}_{M}/2$ be two thresholds to control the number of basis in a single interval.  The parameter $\hat{r}_{M}$ controls the online computational efficiency, and $\hat{r}_{m}$ guides the coarsening to avoid over-refinement. Let $N_{\textrm{iter}}$ be the maximum number of iterations for adaptive refinement. We start from an arbitrary initial partition $\{\mathcal{T}^{(0)}_j\}_{j=1}^{k_0}$. Let the time partition after the $i$-th iteration be $\{\mathcal{T}^{(i)}_j\}_{j=1}^{k_i}$ and the number of bases in $\mathcal{T}_j^{(i)}$ be $r_j^{(i)}$. In the $i$-th iteration, we perform local refinement and coarsening to obtain partition $\{\mathcal{T}_j^{(i)}\}_{j=1}^{k_i}$. We loop over all the intervals from $j=1$ to $j=k_{i-1}$.
\begin{enumerate}
    \item \textbf{Refinement}: If $r^{(i-1)}_j>\hat{r}_{M}$, we uniformly split $\mathcal{T}^{(i-1)}_j$ into two subintervals.
    \item \textbf{Coarsening}: 
    If $r^{(i-1)}_j<\hat{r}_{m}$, we find an integer $0\leq\Delta j\leq k_{i-1}-j$ satisfying
    \begin{equation}
       \Delta j = \argmax_{0\leq \Delta j\leq k_{i-1}-j} \max\{r^{(i-1)}_{j},r^{(i-1)}_{j+1},\dots,r^{(i-1)}_{j+\Delta_j}\}<\hat{r}_{m}. 
    \end{equation}  
    In other words, $\Delta j$ is the maximum number of intervals saisfying that the numbers of bases in all consecutive intervals $\mathcal{T}^{(i-1)}_j,\dots, \mathcal{T}^{(i-1)}_{j+\Delta j}$ are smaller than $\hat{r}_{m}$.
    Then, we merge these intervals according to the following three cases.
    \begin{itemize}
    \item[(a)] When $j+\Delta j+1\leq k_{i-1}$, i.e. $\mathcal{T}_{j+\Delta j}$ is not the last interval, we merge $\mathcal{T}^{(i-1)}_j,\dots,\mathcal{T}^{(i-1)}_{j+\Delta j}$ to $\mathcal{T}^{(i-1)}_{j+\Delta j+1}$ if no refinement conducted in this interval.
    Otherwise, we merge these intervals to the first sub=interval divided from $\mathcal{T}^{(i-1)}_{j+\Delta j+1}$ in the $i$-th iteration.
    \item[(b)] When $j+\Delta j=k_{i-1}$ and $\Delta j>0$, i.e. $\mathcal{T}_{j+\Delta j}$ is the last interval,  we merge $\mathcal{T}^{(i-1)}_j,\dots,\mathcal{T}^{(i-1)}_{j+\Delta j}$.
    \item[(c)] If $j=k_{i-1}$ and $\Delta j=0$, then $\mathcal{T}_{j}^{(i-1)}=\mathcal{T}_{k_{i-1}}^{(i-1)}$ is the last interval, and we merge it to $\mathcal{T}^{(i-1)}_{k_{i-1}-1}$ or the last sub-interval divided from it if refinement is conducted to $\mathcal{T}^{(i-1)}_{k_{i-1}-1}$ in the $i$-th iteration. 
    \end{itemize}
\end{enumerate}
The above algorithm iterates until the numbers of basis functions in all intervals lie between $[\hat{r}_m, \hat{r}_M]$ or the iteration number has reached $N_{\textrm{iter}}$. After refinement and coarsening, we apply POD in each new interval $\mathcal{T}_j$ to generate local bases with the truncation tolerance $\varepsilon_{\textrm{POD}}$.

In cases (a) and (b) of our coarsening strategy, we preferably merge intervals to later intervals based on the intuition that a dynamical system tends to approach its equilibrium over long time if such equilibrium exists. 

The proposed method ensures that the number of bases in each interval lying in a desired range $[\hat{r}_{m},\hat{r}_{M}]$. Furthermore, compared to the adaptive refinement without coarsening in \cite{adpative_time_partition} starting from one initial interval, the coarsening ability of our method allows the use of empirical initial partitioning to reduce iterations for adaptive partitioning by robustly leveraging various initial partitioning.

\textbf{Equilibrium detection.} In long-time simulations, the underlying problem usually goes close to its equilibrium, which often exhibits a low-rank structure. To improve online efficiency, we remove case (c) from our coarsening strategy to allow an extremely small number of bases in the last interval to capture the low-rank structure of the equilibrium. We also correspondingly change the stopping criteria as the number of basis functions except for the last interval lying in the range $[\hat{r}_m,\hat{r}_M]$. 

Long-term simulations can be greatly benefited by this simple modification in terms of time and memory efficiencies, as it allows using a very low-dimensional reduced order space when the underlying problem approaches its equilibrium.

\subsection{Hybrid ROM with time partitioning}\label{sec:hybridPOD} 
In extreme scenarios, low-rank structures of the solution manifold may not be effectively captured by piecewise linear in time approximations. POD with adaptive time partitioning may require excessive refinement to maintain accuracy. For example, the snapshot matrix for the propagation of a delta function could be an identity matrix.  

To address this issue, autoencoder-based ROM \cite{lee2020model,choi2021space} is an attractive alternative to POD, due to its potential to capture the solution dynamics with a latent space whose dimension is as low as the intrinsic dimension of the solution manifold \cite{JSC}. However, training an autoencoder is much more expensive than computing SVD. To avoid over-refinement without dramatically compromising the offline efficiency, we only employ the autoencoder when necessary. Additionally, this strategy of applying autoencoder locally is able to reduce the variance and amount of data, which facilitates an easier training process.

\textbf{Goal-oriented hybridization.}
Our hybrid method still follows a goal-oriented strategy, where the complexity of the algorithm is controlled by the choices of parameters to meet the specific time and memory requirement in practice. Define $\tau_{\min}\geq 0$ as the minimal size of time intervals allowed to avoid over-refinement. The maximum number of basis functions allowed in each interval is denoted as $\hat{r}_{M}$.
We will apply the autoencoder in a time interval $\mathcal{T}_j$, if the following two conditions are satisfied simultaneously.
\begin{enumerate}
    \item The number of reduced basis functions required in $\mathcal{T}_j$ is larger than $\hat{r}_{M}$. 
    \item Uniform partitioning $I_j$ will result in subintervals whose sizes are smaller than $\tau_{\min}$.
\end{enumerate}

\textbf{Structure of convolutional autoencoder. }
In what follows, we conclude with the architecture of the convolutional autoencoder.
Autoencoder is a neural network providing a low-rank approximation to the identity mapping 
\begin{align*}
\Id:\mathbb{R}^N&\to\mathbb{R}^N,\\
    \textbf{y}&\mapsto \textbf{y}.
\end{align*}
An autoencoder consists of two parts: an encoder neural network $\mathcal{E}_\theta$ and a decoder neural network $\mathcal{D}_\theta$, where $\theta$ represents trainable parameters. The encoder defines a nonlinear mapping compressing a high-dimensional input $y\in\mathbb{R}^N$ to a low-dimensional reduced order representation $y_r\in\mathbb{R}^r$ with $r\ll N$. The decoder is a nonlinear mapping that reconstructs the high-dimensional input $y\in\mathbb{R}^N$ from its corresponding low-dimensional latent representation $y_r\in\mathbb{R}^r$. In other words, $y\approx\mathcal{D}_{\theta}(y_r)=\mathcal{D}_{\theta}\circ\mathcal{E}_{\theta}(y)$. 

\begin{figure}
    \centering
    \includegraphics[width=1\linewidth]{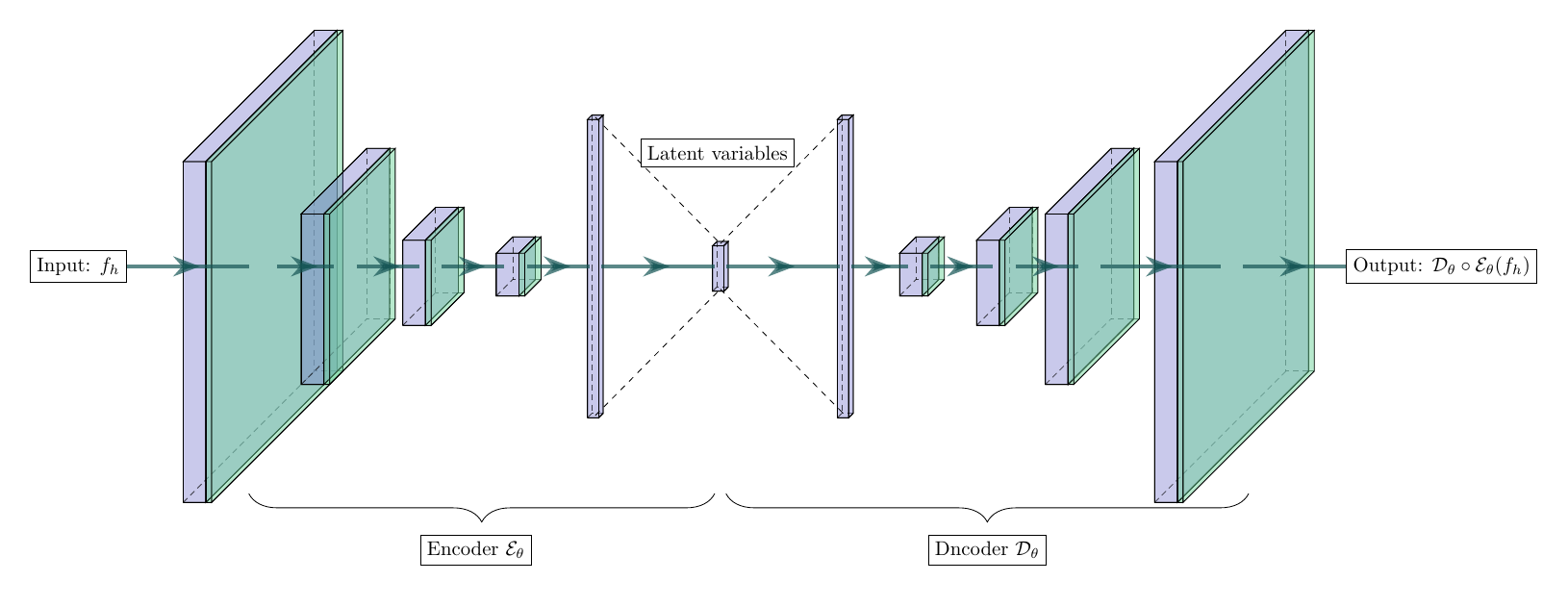}
    \caption{Architecture of convolutional autoencoder: the purple blocks represent the convolutional and transposed convolutional layers, while the green blocks represent the activation function.}
    \label{fig:AE_arch}
\end{figure}

Our training data, i.e. high-fidelity solutions, are nodal values determined by the numerical solution given by our discontinuous Galerkin solver on rectangular structured meshes. Due to the structured nature of our data, Following \cite{lee2020model,CAE_LSTM,DUAN2024112621}, we employ a convolutional autoencoder where feature dimensions correspond to mesh resolution, with $n_v$ channels corresponding to each velocity (or 1 channel for $\rho$).
Key advantages of adopting convolutional architecture are as follows.
\begin{enumerate}
    \item Weight-sharing dramatically reduces the number of trainable parameters, which facilitates better efficiency compared to fully-connected networks.
    \item Local kernel operations capture neighborhood correlations, mirroring the stencil structure of our numerical solution.
\end{enumerate}

The neural network architecture is illustrated in Figure \ref{fig:AE_arch}. 
Specifically, the encoder $\mathcal{E}_\theta$ comprises a sequence of convolutional layers followed by a fully-connected layer. The decoder $\mathcal{D}_\theta$ 
mirrors this architecture, beginning with a fully-connected layer and followed by a sequence of transposed convolutional layers. In each layer, we use hyperbolic tangent (Tanh) as the nonlinear activation function.

Given the snapshots $\{f_h(t_i;\mu_j): i=1,\dots,n_t; j=1,\dots, n_p\}$ whose elements correspond to nodal values in the physical space, we determine the trainable parameters $\theta$ via minimizing the following loss function:
\begin{align}\label{eq:loss}
    \mathcal{L} = \frac{1}{n_t n_p}\sum_{i=1}^{n_t}\sum_{j=1}^{n_p} \|f_h(t_i, \mu_j)-\mathcal{D}_\theta\circ\mathcal{E}_\theta(f_h(t_i, \mu_j))\|^2_{l^2}
\end{align}
The neural network is constructed and implemented through {\tt{Pytorch}} library and trained by minimizing the loss function \eqref{eq:loss} with an ADAM optimizer \cite{kingma2014adam}.

\section{Numerical experiments}\label{sec:num_experiments}
In this section, we present numerical tests of the proposed model reduction methods on three examples.
The following relative errors will be evaluated in various settings at a given testing time and parameter $(t_{\text{test}}, \mu_{\text{test}})$
\begin{align*}
    e_f(t_{\text{test}},\mu_{\text{test}}) &=\frac{\|f_h(t_{\text{test}},\mu_{\text{test}})-\tilde{f}_h(t_{\text{test}},\mu_{\text{test}})\|_{l^2}}{\|f_h(t_{\text{test}},\mu_{\text{test}})\|_{l^2}}, \\
    e_{\rho}(t_{\text{test}},\mu_{\text{test}}) &=\frac{\|\rho_h(t_{\text{test}},\mu_{\text{test}})-\tilde{\rho}_h(t_{\text{test}},\mu_{\text{test}})\|_{l^2}}{\|\rho_h(t_{\text{test}},\mu_{\text{test}})\|_{l^2}},
\end{align*}
where $f_h$ and $\rho_h$ are high-fidelity full order solutions, $\tilde{f}_h$ and $\tilde{\rho}_h$ are reduced order solutions. Note that in our implementation, $t_{\text{test}}$ falls within training sampling times, while $\mu_{\text{test}}$ lies outside training sampling parameters. The average errors of the whole testing set are defined as
\begin{align*}
    E_f = \frac{1}{N_{t_{\text{test}}}N_{\mu_{\text{test}}}}\sum_{t_{\text{test}}}\sum_{\mu_{\text{test}}}e_f(t_{\text{test}}, \mu_{\text{test}}),\\
    E_\rho = \frac{1}{N_{t_{\text{test}}}N_{\mu_{\text{test}}}}\sum_{t_{\text{test}}}\sum_{\mu_{\text{test}}}e_{\rho}(t_{\text{test}}, \mu_{\text{test}}),
\end{align*}
where $N_{t_{\text{test}}}$ and $N_{\mu_{\text{test}}}$ are the numbers of the testing times and parameters, respectively. 

In our experiments, we choose the reference solution $\tilde{f}_h$ to be the mean of all $n_s$ snapshots. We emphasize that data normalization is essential for the efficiency of the neural network training phase. Previous studies, such as \cite{lee2020model,JSC}, typically rescale snapshots to the range $[0,1]$ using the same maximum and minimum values for all elements of the snapshot matrix. However, due to the multiscale nature of the RTE and the variability introduced by different sampling times and parameters, our training snapshots exhibit significant variance, making this rescaling approach less effective. Therefore, we adopt an alternative normalization approach that rescales the input snapshots to the range $[-1,1]$, aligning with the range of the hyperbolic tangent activation function. We first reshape the snapshot matrix $S\in\mathbb{R}^{(n_\mathbf{x}n_\mathbf{v})\times (n_\mu n_t)}$ into $W\in\mathbb{R}^{(n_\mathbf{x}n_\mu)\times(n_\mathbf{v}n_t)}$. Then we define
\begin{align*}
    W_{\max}^j = \max_{i=1,\dots,n_\mathbf{x}n_\mu} W_{ij} \in\mathbb{R},\quad W_{\min}^j = \min_{i=1,\dots,n_\mathbf{x}n_\mu} W_{ij} \in\mathbb{R}.
\end{align*}
so that the following transformation is applied to reshaped snapshot matrix $W$
\begin{align}\label{eq:data_normalize}
    W_{ij}\mapsto2\frac{W_{ij}-W_{\min}^j}{W^j_{\max}-W^j_{\min}}-1,\quad i =1,2,\dots,n_\mathbf{x}n_\mu, \quad j = 1,2,\dots, n_\mu n_t.
\end{align}
The same transformation is applied to the testing snapshots as well, but the corresponding maximum and minimum values are still computed over training data.

The ground truth snapshots are generated using upwind linear DG method with backward Euler time discretization. All the numerical experiments are implemented through Python. Example $1$ is tested on an iMac with M1 chip. For example $2$, we train the autoencoder on the \verb|Tianhe-2| supercomputer with $1$ GPU core and $6$ CPU cores, while all other computations were performed on a MacBook Air equipped with an Apple M1 chip and $8$ CPU cores. All experiments for Example $3$ were performed on a system equipped with an NVIDIA GeForce RTX 4090 GPU and a 16-core CPU.

\subsection{Example 1: two-material problem}
\begin{figure}[htbp]
\centering
\includegraphics[width=0.75\textwidth]{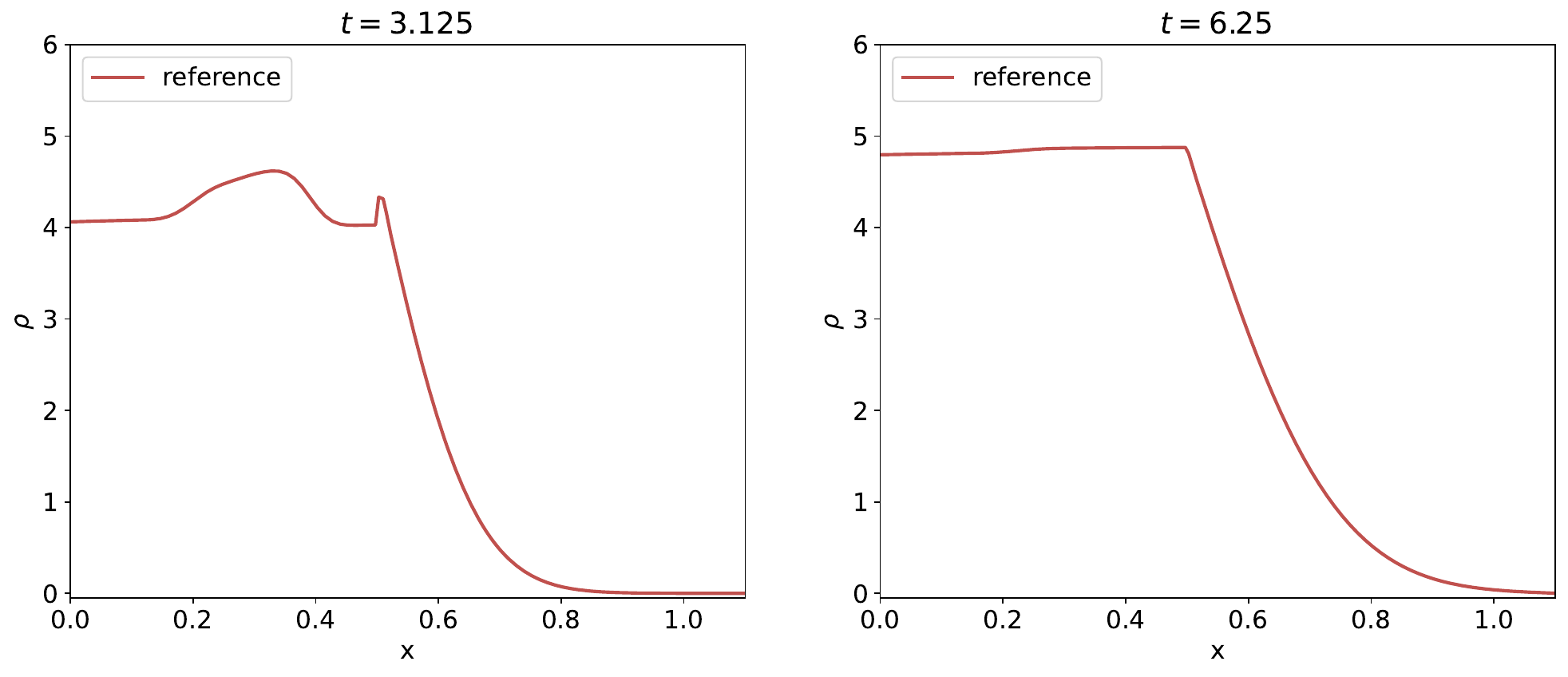}
\includegraphics[width=0.75\textwidth]{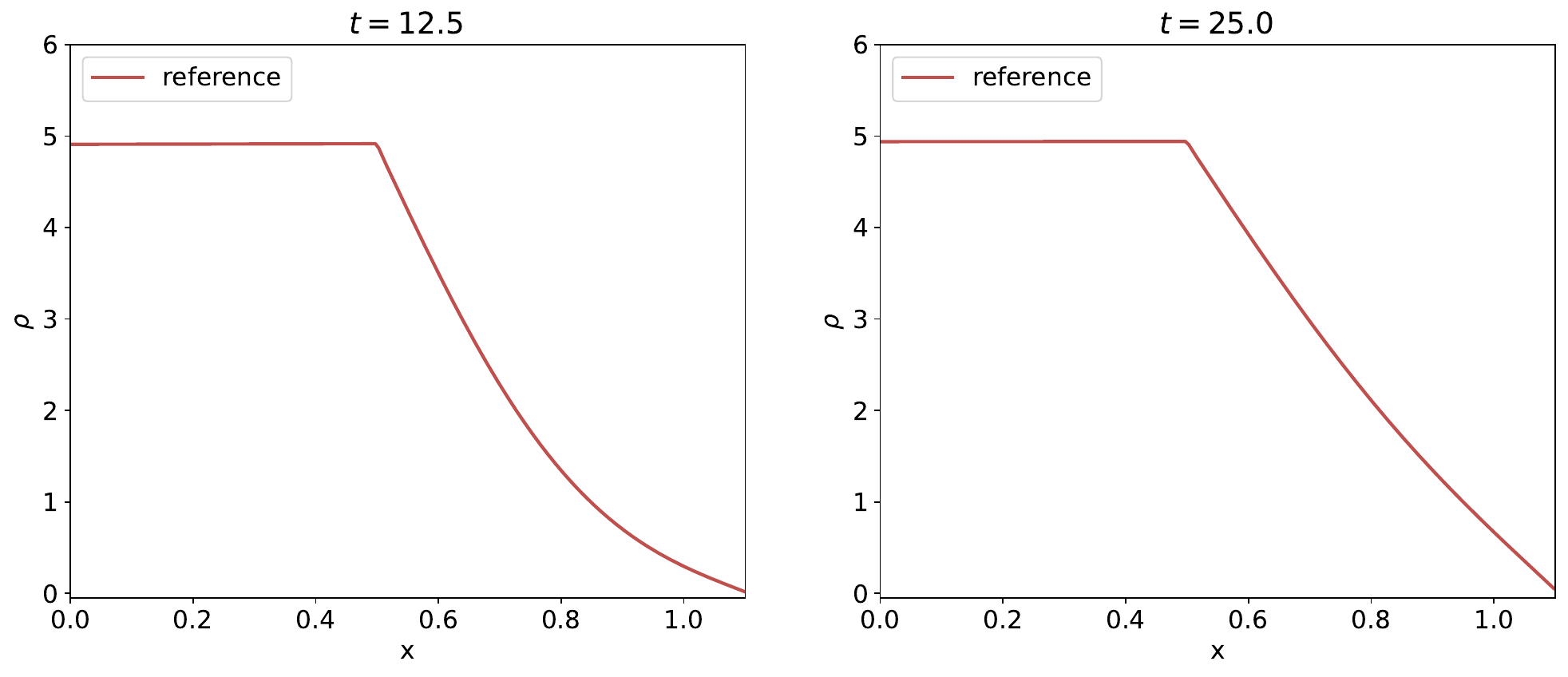}
	\caption{Reference solutions for Example 1 corresponding to $\mu=5$.}
	\label{fig:eg1_reference_sol}
\end{figure}

We first consider a  two-material problem on $\Omega_x=[0,1.1]$ with inflow boundary condition and zoer initial condition in the 1D slab geometry: 
\begin{align*}
\sigma_a(x)=0, \quad\sigma_s(x)=\begin{cases}
0,\hspace{11pt}\quad\quad 0\leq x\leq1,\\
100, \quad\quad 1< x\leq1.1,
\end{cases}\\
f(0,v,t)=\mu \text{ for }v>0; \quad f(1.1,v,t)=0 \text{ for }v<0.
\end{align*}
Here, parameter $\mu\in[4,6]$ determines the left inflow boundary conditions. At the discrete level, we uniformly partition the computational domain into $88$ elements, using $16$ Gauss-Legendre points in the angular space and solve the problem from $t=0$ to $t=25$ with time step size $\Delta t=\frac{1}{80}$. Reference solutions for $\mu=5$ at four sampling times are given in Figure \ref{fig:eg1_reference_sol}.

We uniformly sample $11$ training parameters $\mu_j=4+0.2j$ with $j=0,\dots,10$, and build our reduced order models using all snapshots corresponding to these training parameters and $0\leq t \leq 25$. We set the truncation tolerance in POD method as $\varepsilon=10^{-4}$, and the threshold for number of basis of each interval in adaptive partitioning as  $\hat{r}_{M}=15$ and $\hat{r}_{m}=5$. In the online prediction, we test the performance of the proposed method with $5$ random test parameters uniformly sampled from $[4,6]$.

\begin{figure}[htbp]
\centering
\includegraphics[width=0.49\textwidth]{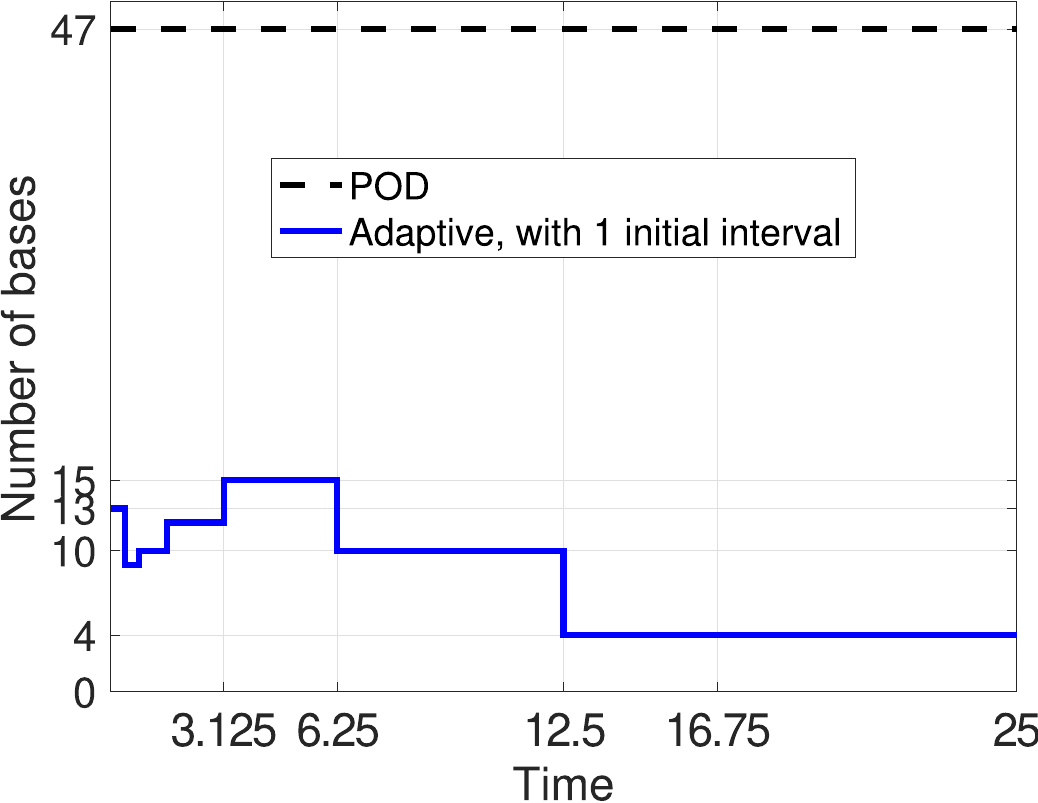}
\includegraphics[width=0.49\textwidth]{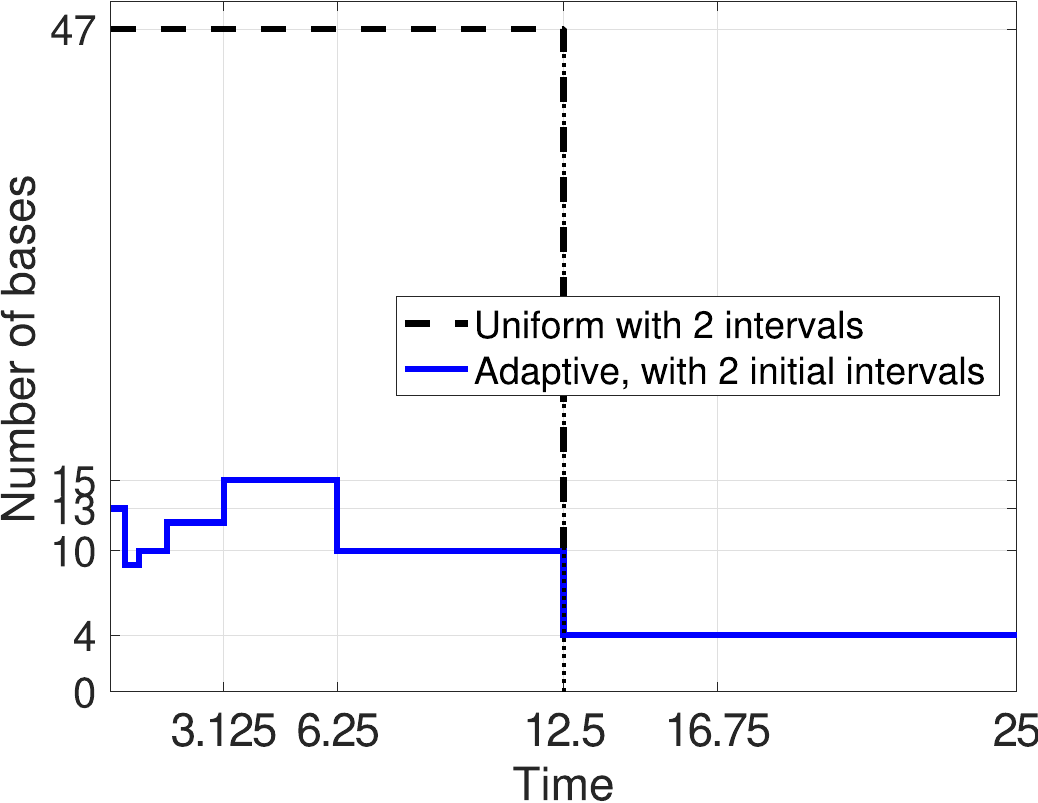}
\includegraphics[width=0.49\textwidth]{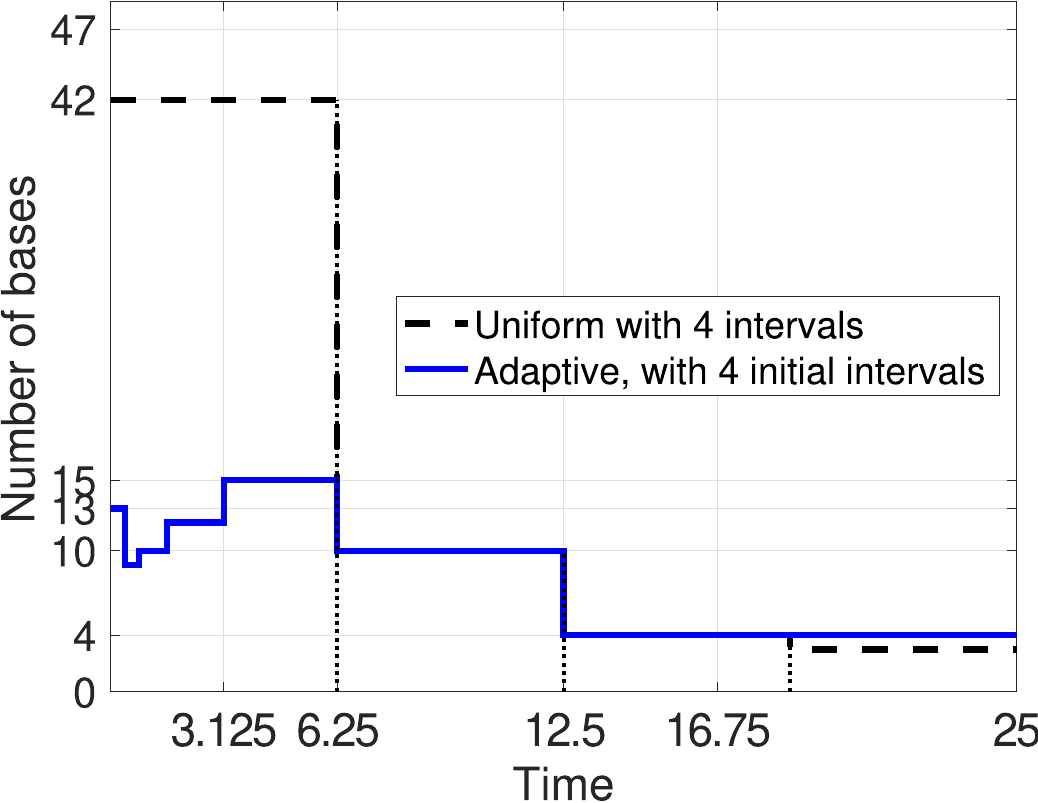}
\includegraphics[width=0.49\textwidth]{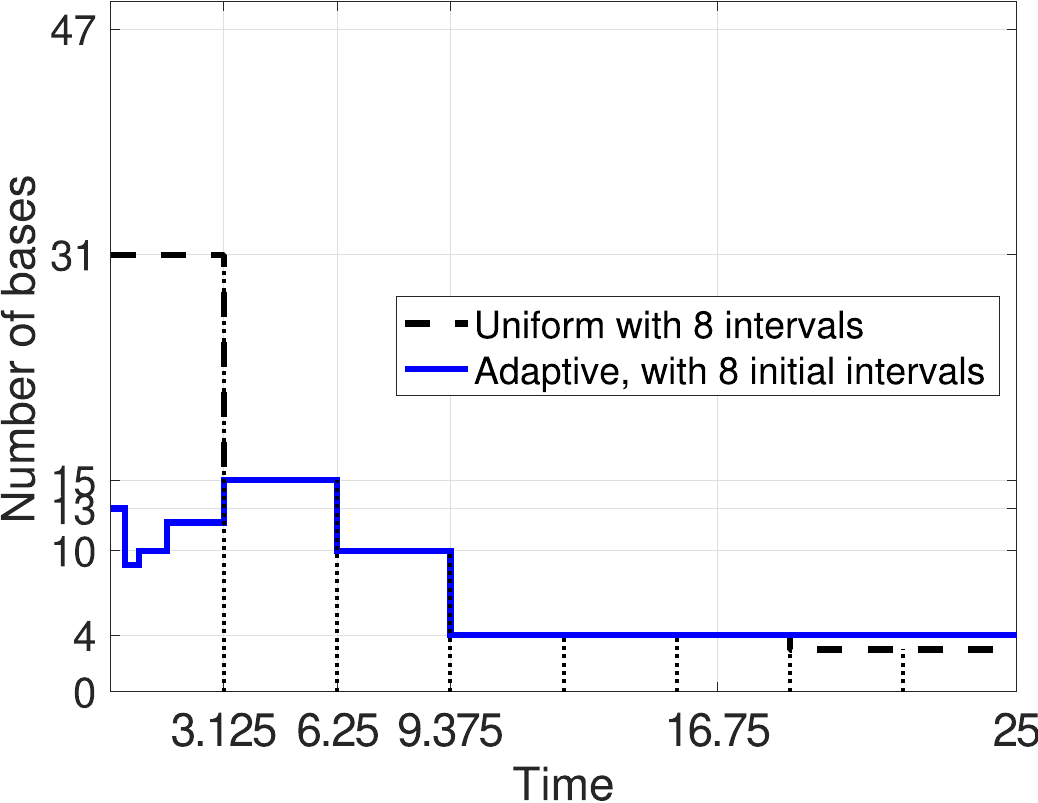}
	\caption{Number of basis for various adaptive time partitioning strategies for Example 1.}
	\label{fig:eg1_partition}
\end{figure}

\textbf{Robustness of adaptive partitioning.} To demonstrate the robustness of adaptive partitioning, we test its performance with various uniform initial partitioning using $1$, $2$, $4$ and $8$ subintervals.
The resulting time partition and the number of bases in each interval are presented in Figure \ref{fig:eg1_partition}. 

We observe that, regardless of the initial time partition, our adaptive approach uses more refined intervals at first and gradually switches to largely coarsened intervals in the evolution of the system. This phenomenon aligns well with the behavior of the system at hand. Initially, there is no particle in the computational domain. As particles enter the domain from the left boundary, they first go through a free-streaming region. Then, an internal layer is formed when particles start to enter the strong scattering region on the right side. As time passes, the particle system gradually approaches its steady state. In the end, the rank of the system decreases and does not change significantly over time.  

\begin{figure}[htbp]
\centering
\adjustbox{valign=t}{\includegraphics[width=0.49\textwidth]{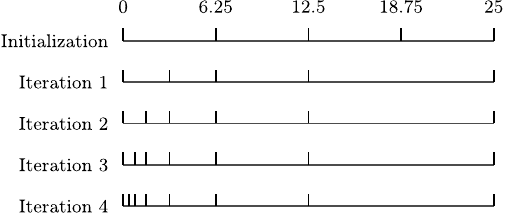}}
\raisebox{-19mm}{\includegraphics[width=0.49\textwidth]{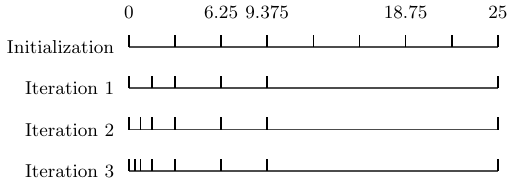}}
	\caption{How time partition evolves in Example 1 with different uniform initial partitions. Left: $4$ initial intervals. Right: $8$ initial intervals.}
	\label{fig:eg1_iter}
\end{figure}

Another interesting observation is that the final partitions of our adaptive approach remain consistent when the initial uniform partition has $1$, $2$ or $4$ intervals. The difference between them and the final partition starting from $8$ uniform intervals is that the last coarsened long interval starts earlier when $8$ intervals. As shown in Figure \ref{fig:eg1_iter}, the first iteration with $8$ intervals is able to coarsen the last interval more aggressively due to its better resolution in the initial partition. 

Overall, the proposed adaptive method is capable of partitioning the time domain robustly and adaptively with various initial partition. 
 
\textbf{Comparison between uniform and adaptive partitioning.}
We first test the accuracy of uniform and adaptive partitioning methods. As shown in Table \ref{table:eg1_error}, the average prediction errors for testing parameters are comparable across various methods. 

Table \ref{table:eg1_acceleration} presents the relative average numbers of reduced basis functions required for different methods with respect to POD (i.e. $\frac{\text{average number of bases}}{\text{number of basis required for POD}}$) in each interval, while Table \ref{table:eg1_online_time} shows the online computational time. 
As shown in Table \ref{table:eg1_acceleration}, Table \ref{table:eg1_online_time} and Figure \ref{fig:eg1_partition}, uniform partitioning strategy under-refines the crucial initial period $t\in[0,3.125]$, and hence only leads to marginal acceleration compared to POD in $t\in[0,3.125]$. In contrast, our adaptive approach sufficiently refines this time interval and significantly reduces the average number of basis functions in this period to less than $25\%$ compared to POD, resulting in close to $3$ times acceleration over POD in the corresponding region. We suspect that the mismatch between the reduction of the number of basis and the acceleration is due to the overhead to transit from one interval to another in the adaptive method. 

Uniform time partition with $8$ intervals not only under-refines $t\in[0,3.125]$, but also unnecessarily over-refines $t\in[12.5,25]$. We highlight that the adaptive approach successfully avoids both issues.

In summary, our adaptive method is able to avoid both the under-refinement during the initial transient dynamics and the over-refinement in long-time as the system approaches equilibrium.

\begin{table}[htbp]
\centering
\footnotesize
\begin{tabular}{|c|c|c|c|c|c|c|}
\hline & $E_f$ & $E_\rho$ \\ \hline 
POD         & $2.35\times 10^{-5}$ & $1.15\times10^{-5}$ \\ \hline
Uniform-$2$ & $2.04\times 10^{-5}$ & $1.01\times10^{-5}$\\ \hline
Uniform-$4$ &  $1.98\times 10^{-5}$ & $1.15\times 10^{-5}$\\ \hline
Uniform-$8$ &  $2.07\times 10^{-5}$ & $1.53\times 10^{-5}$\\ \hline \midrule \hline
Adaptive-$4$ & $1.37\times10^{-5}$ & $6.05\times10^{-6}$ \\ \hline
Adaptive-$8$ &  $1.15\times10^{-5}$ & $8.20\times10^{-6}$ \\ \hline
\end{tabular}
\caption{Average relative errors for the testing parameters for Example 1. Uniform-$\theta$: uniform partitions with $\theta$ intervals. Adaptive-$\theta$: adaptive partitioning with $\theta$ uniformly distributed initial intervals.}
\label{table:eg1_error}
\end{table}

\begin{table}[htbp]
\centering
\begin{tabular}{|c|c|c|c|c|c|c|}
\hline
   & $t\in[0,3.125]$ & $t\in[3.125,6.25]$ & $t\in[6.25,12.5]$ & $t\in[12.5,25]$                                              \\ \hline
POD & $1.00$ & $1.00$ & $1.00$ & $1.00$  \\ \hline
Uniform-$4$ & $0.89$ & $0.89$ & $0.21$& $0.07$ \\ \hline
Uniform-$8$ & $0.66$ & $0.32$ & $0.15$ & $0.07$ \\ \hline \midrule \hline
Adaptive-$4$ & $0.24$ & $0.32$ & $0.21$ & $0.09$ \\ \hline
Adaptive-$8$ &  $0.24$ & $0.32$ & $0.15$ & $0.09$  \\ \hline
\end{tabular}
\caption{Relative average number of basis with respect to POD (i.e. $\frac{\text{average number of bases}}{\text{number of basis required for POD}}$) in Example 1. Uniform-$\theta$: uniform partitions with $\theta$ intervals. Adaptive-$\theta$: adaptive partitioning with $\theta$ uniformly distributed initial intervals.}
\label{table:eg1_acceleration}
\end{table}
\begin{table}[htbp]
\centering
\begin{tabular}{|c|c|c|c|c|c|c|}
\hline
   & $t\in[0,3.125]$ & $t\in[3.125,6.25]$ & $t\in[6.25,12.5]$ & $t\in[12.5,25]$ & $t\in[0,25]$                                               \\ \hline
POD         & $0.55$ & $0.43$ & $0.89$ & $1.68$ & $3.55$  \\ \hline
Uniform-$2$ & $0.55$ & $0.43$ & $0.89$ & $0.22$ & $2.09$  \\ \hline
Uniform-$4$ & $0.50$ & $0.41$ & $0.29$& $0.22$  & $1.38$\\ \hline
Uniform-$8$ & $0.39$ & $0.23$ & $0.25$ & $0.27$ & $1.14$ \\ \hline \midrule \hline
Adaptive-$4$ & $0.18$ & $0.21$ & $0.29$ & $0.29$ & $0.99$ \\ \hline
Adaptive-$8$ & $0.19$ & $0.24$ & $0.23$ & $0.26$   & $0.91$\\ \hline
\end{tabular}
\caption{Online computational time (sec) in Example 1. Uniform-$\theta$: uniform partitions with $\theta$ intervals. Adaptive-$\theta$: adaptive partitioning with $\theta$ uniformly distributed initial intervals.}
\label{table:eg1_online_time}
\end{table}
\textbf{Strength and weakness of the hybridized strategy.}
We then compare the performance of the hybrid method and POD with adaptive and uniform time partitioning. 

Here, we set the truncation tolerance in POD as $\varepsilon=5\times10^{-3}$  and the threshold for the number of bases of each interval in adaptive partitioning as $\hat{r}_{M}=15$ and $\hat{r}_{m}=5$. The minimum time interval allowed in the hybrid method is set as $\tau_{\min}=6.25$.
We use $4$ intervals in the uniform partitioning method, and employ them as the initial partition for the adaptive and hybrid approach.

In the hybrid method, we set the dimension of the latent space constructed by the autoencoder as $r=8$. We refer to Appendix \ref{apx:ae} for a detailed setup of hyperparameters and training of our autoencoder. 

\begin{figure}[htbp]
\centering
\begin{tikzpicture}
    \node[anchor=south west,inner sep=0] (image) at (0,0) {\includegraphics[width=0.49\textwidth]{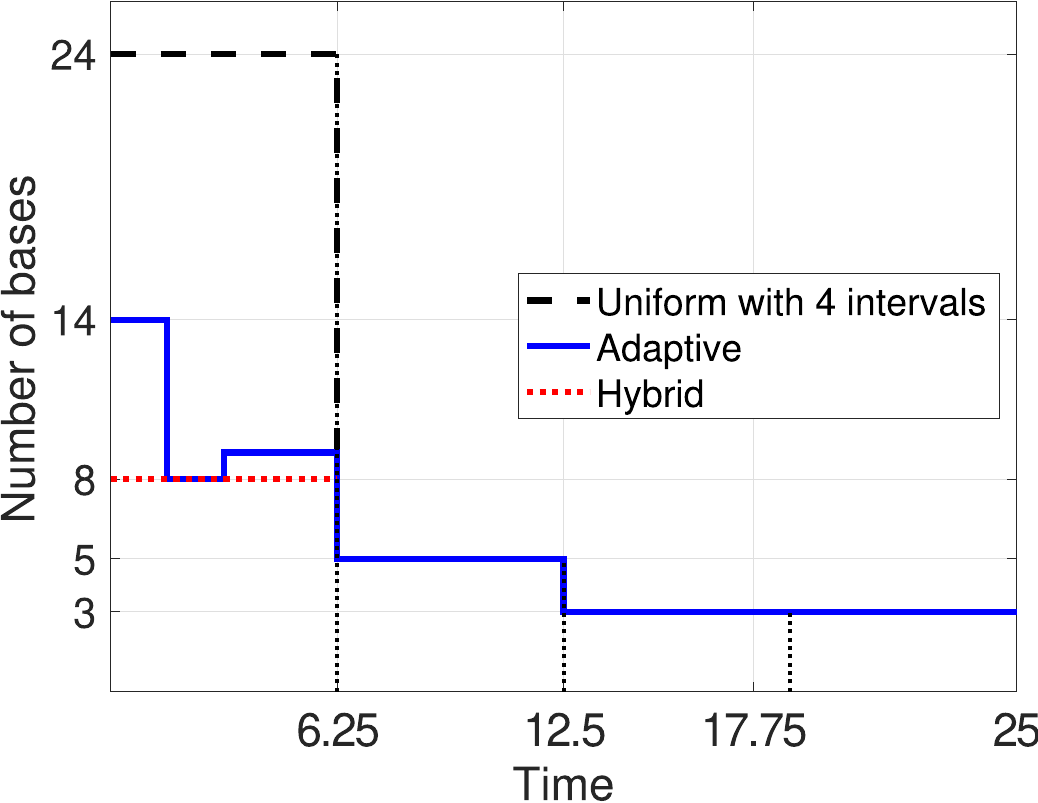}};
    
    \node[red,font=\footnotesize] at (1.5,1.7) {AE};
    \node[red,font=\footnotesize] at (1.5,1.3) {applied};
    
\end{tikzpicture}
\caption{Dimensions of the latent spaces of various methods in Example 1, with the POD truncation tolerance $\varepsilon_{\text{POD}}=5\times10^{-3}$. AE: autoencoder.}
\label{fig:eg1_hybrid_dim}
\end{figure}

The dimensions of the latent spaces for various methods are presented in Figure \ref{fig:eg1_hybrid_dim}.
The hybrid method results in a three-interval partition with $\mathcal{T}_1=[0,6.25]$, $\mathcal{T}_2=[6.25,12.5]$ and $\mathcal{T}_3=[12.5,25]$, while the adaptive POD leads to a five-interval partition. The autoencoder is applied only in the first time interval ($[0,6.25]$) in the hybrid approach. As shown in Table \ref{table:eg1_hybrid_time}, by introducing the autoencoder, the hybrid ROM results in approximately $1.82$ and $3.24$ times speed-up versus POD with adaptive and uniform partitioning. Meanwhile, Table \ref{table:eg1_hybrid_error} shows the relative errors in the time interval $t\in[0,6.25]$ where the autoencoder is applied in the hybrid method. As we can see from the table, the three approaches achieve comparable accuracy.

The hybrid ROM leverages an autoencoder to achieve more effective solution compression without extensive time partitioning. However, it also presents a trade-off: When training the autoencoder, a highly nonconvex optimization problem needs to be solved. Compared to SVD, solving this optimization problem is more time-consuming. Also, the parameters of the trained autoencoder are suboptimal with relatively limited accuracy due to the highly nonconvexity of the neural network, while SVD provides the theoretically guaranteed optimal linear approximation. 

\begin{table}[htbp]
\centering
\begin{tabular}{|c|c|c|c|c|c|c|}
\hline
  & $t\in[0,6.25]$ & $t\in[6.25,12.5]$ & $t\in[12.5,25]$ & $t\in[0,25]$                                               \\ \hline
Uniform  & $0.55$ & $0.18$ & $0.18$ & $0.90$ \\ \hline
Adaptive & $0.31$ & $0.17$ & $0.19$ & $0.67$ \\ \hline 
Hybrid & $\mathbf{0.17}$ (AE applied) & $0.17$ & $0.19$ & $\mathbf{0.53}$ \\ \hline
\end{tabular}
\caption{Comparison of online computational time (sec) between the hybrid ROM and POD with uniform and adaptive partitioning in Example 1. Truncation tolerance of POD is set as $\varepsilon_{\text{POD}}=5\times10^{-3}$.}
\label{table:eg1_hybrid_time}
\end{table}
\begin{table}
\centering
\begin{tabular}{|c|c|c|c|c|c|c|}
\hline
  & $E_f$ for $t\in[0,6.25]$ &  $E_\rho$ for $t\in[0,6.25]$                                      \\ \hline
Uniform  & $2.19\times10^{-3}$ & $9.42\times10^{-4}$   \\ \hline
Adaptive & $1.78\times10^{-3}$ & $7.04\times10^{-4}$  \\ \hline 
Hybrid & $2.56\times10^{-3}$ & $1.21\times10^{-3}$  \\ \hline
\end{tabular}
\caption{The average relative prediction error for the hybrid ROM and POD with uniform and adaptive partitioning in Example 1. Truncation tolerance of POD is set as $\varepsilon_{\text{POD}}=5\times10^{-3}$.}\label{table:eg1_hybrid_error}
\end{table}


\subsection{Example 2: three-material problem}
In the second example, we consider a $1$-dimensional problem on $\Omega_x=[0,2]$ and $\Omega_v=[-1,1]$ with a source term. The setup is as follows.
\begin{align*}
    \sigma_a(x)&=0, \quad\sigma_s(x)=\begin{cases}
0,\quad\quad |x-1|\leq 0.3,\\
\mu, \quad\quad 1.3<x\leq2,\\
5, \quad\quad 0\leq x<0.7,
\end{cases}\\
G(x,v,t) &= \begin{cases}
1,\quad\quad |x+1|\leq 0.5,\\
0, \quad\quad |x+1|>0.5,\\
\end{cases}\\
f(0,v,t)&=0, \quad f(2,v,t)=0; \quad f(x,v,0)=10^3 e^{-\frac{(x-1)^2}{10^{-6}}}.
\end{align*}
We solve the equation up to the final time $T=20$. 

We use a uniform mesh with $128$ elements to partition the computational domain $\Omega_x=[0,2]$, and $16$ Gauss-Legendre points to discretize the angular space $\Omega_v=[-1,1]$. Therefore, in this example, the numbers of spatial and angular degrees of freedom are $n_x = 2\times128$ and $n_v = 16$. The time step size is set to $\Delta t=0.01$. We generate training snapshots at $n_t=200$ uniform sampling times in $[0.01,20]$ with $n_p=64$ parameter values $\mu=\{80,81,\dots,143\}$. We test the performance of our proposed methods against the snapshots with $10$ randomly sampled $\mu\in[75,150]$ from a uniform distribution. In this example, the cubic piecewise polynomial interpolation scheme is adopted in the online extrapolation.

\begin{figure}[htbp]
	\centering
   	\includegraphics[width=1\textwidth]{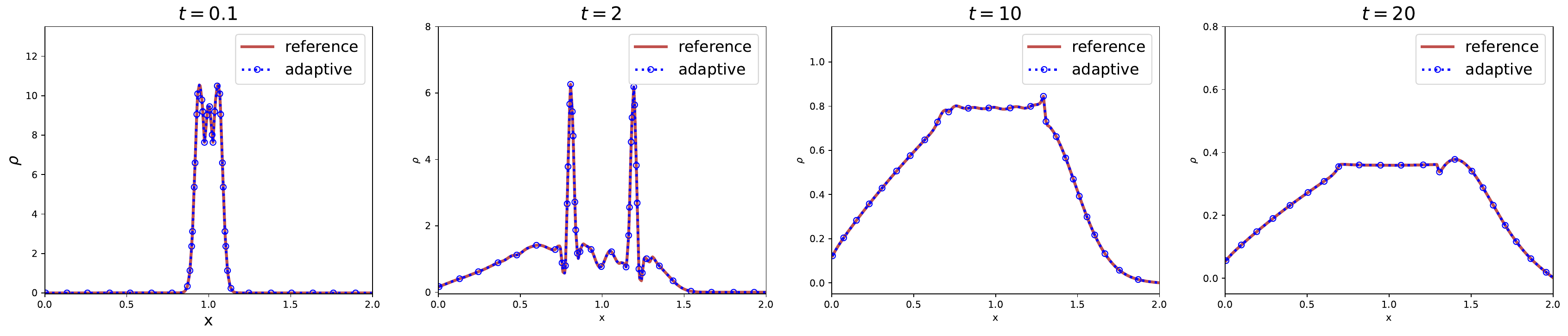}
	\caption{$\rho$ in Example 2 reconstructed by adaptive method with $\hat{r}_{M}=20$ and $\hat{r}_{m}=10$, $x\in[0,2]$, $\mu=136.10$, $\varepsilon_{\text{POD}}=10^{-4}$.}
	\label{fig:eg2_rho}
\end{figure}

\textbf{Comparison between uniform and adaptive partitioning.}
We start by studying the performance of uniform and adaptive time partitioning POD. 
classical POD and uniform partitioning into $k=10$ intervals are chosen as the baseline methods, while the adaptive partitioning method with $\hat{r}_M = 20$ and $\hat{r}_m$ is examined against these baseline methods. We show in Table \ref{table:eg2_uni_adap} the relative error $E_f$, $E_{\rho}$ and online CPU time against $10$ testing parameters in $[75,150]$. It can be witnessed that our adaptive method achieves the best accuracy and the most efficient online stage compared to the baseline methods. The reconstructed $\rho$ using the adaptive method on the test parameter $\mu=136.10$ at $t=0.1,2,10,20$ is shown in Figure \ref{fig:eg2_rho}, where all predicted solutions align well with the reference solutions.
\begin{table}[htbp]
\centering
\begin{tabular}{|c|c|c|c|}
\hline
         & $E_f$    & $E_\rho$ & online time (sec) \\ \hline
POD      & $1.88\times10^{-4}$ & $8.24\times10^{-5}$ & $5.96$             \\ \hline
Uniform-$10$  & $7.77\times10^{-5}$ & $4.91\times10^{-5}$ & 1.38             \\ \hline
Adaptive & \textbf{$6.46\times10^{-5}$} & \textbf{$4.04\times10^{-5}$} & \textbf{$1.18$}\\ \hline
\end{tabular}
\caption{Comparison of relative errors and online time of POD, uniform partitioning and adaptive partitioning methods in Example 2.}
\label{table:eg2_uni_adap}
\end{table}

\begin{figure}[htbp]
\centering
	\includegraphics[width=0.5\textwidth]{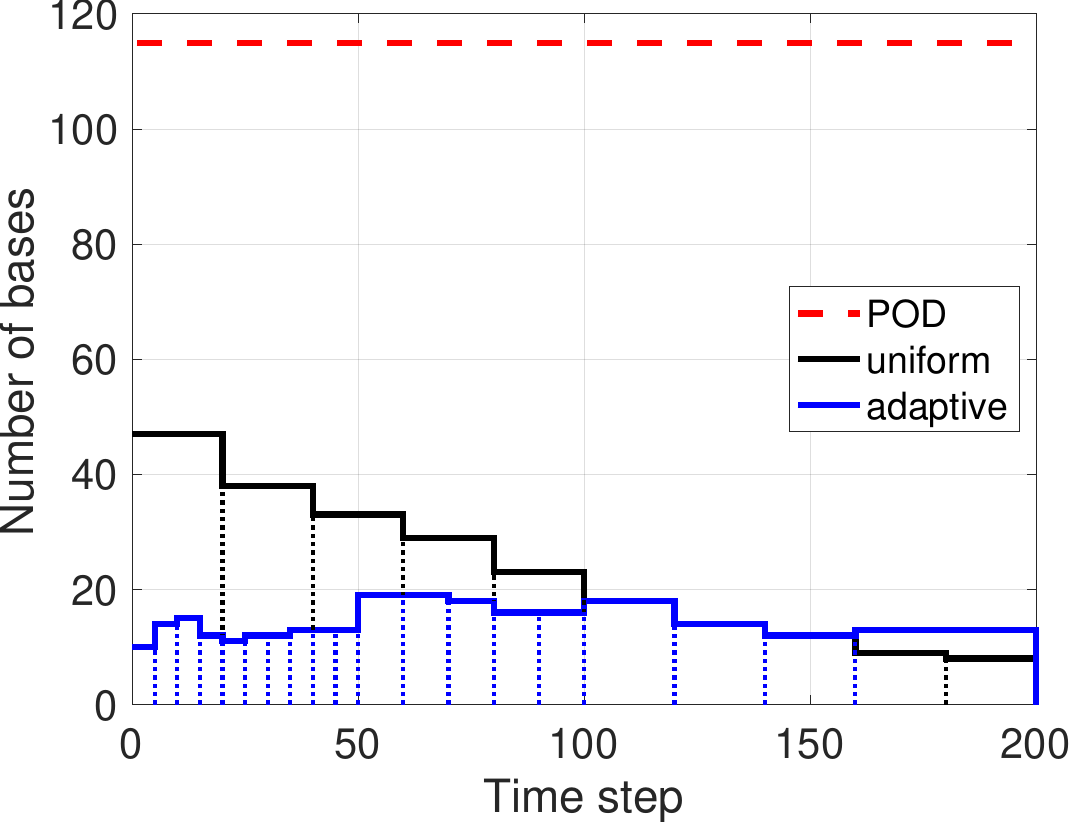}
	\caption{Illustration of number of basis functions on each interval when using different methods in Example 2.}
	\label{fig:eg2_dim}
\end{figure}

Figure \ref{fig:eg2_dim} displays the number of bases required in each time interval using different methods. It shows that the adaptive approach effectively regulates the maximum and minimum number of bases via $\hat{r}_{M}$ and $\hat{r}_{m}$, respectively, ensuring the ability of avoiding both under- and over-refinements. As we can see from the figure, uniform piecewise linear in time partitioning proves inadequate in the initial time intervals, and the adaptive method compensates through finer time partition; conversely, as the system approaches equilibrium, the adaptive method employs coarser partitions to optimize memory usage without compromising accuracy in these regions.

\textbf{Strength and weakness of the hybridized strategy.}
We then study the online efficiency of the autoencoder hybrid method, as described in Section \ref{sec:hybridPOD}, in the current example with input parameters $\tau_{\min}=10, \hat{r}_{M}=15$.

\begin{figure}[htbp]
	\centering
       \subfigure[]{
		\begin{minipage}{1\textwidth}
			\includegraphics[width=1\textwidth]{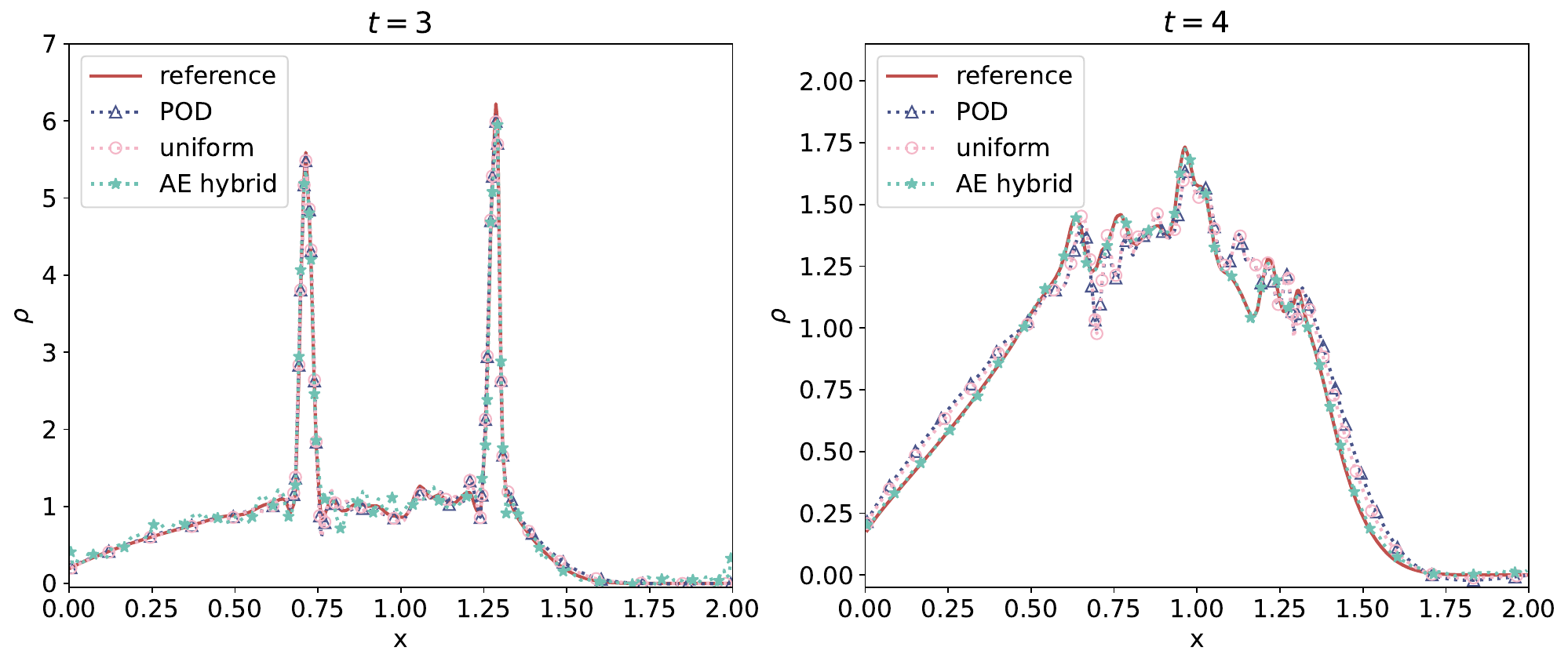}
		\end{minipage}
	}%
    
    	\subfigure[]{
    		\begin{minipage}{1\textwidth}
   		 	\includegraphics[width=1\textwidth]{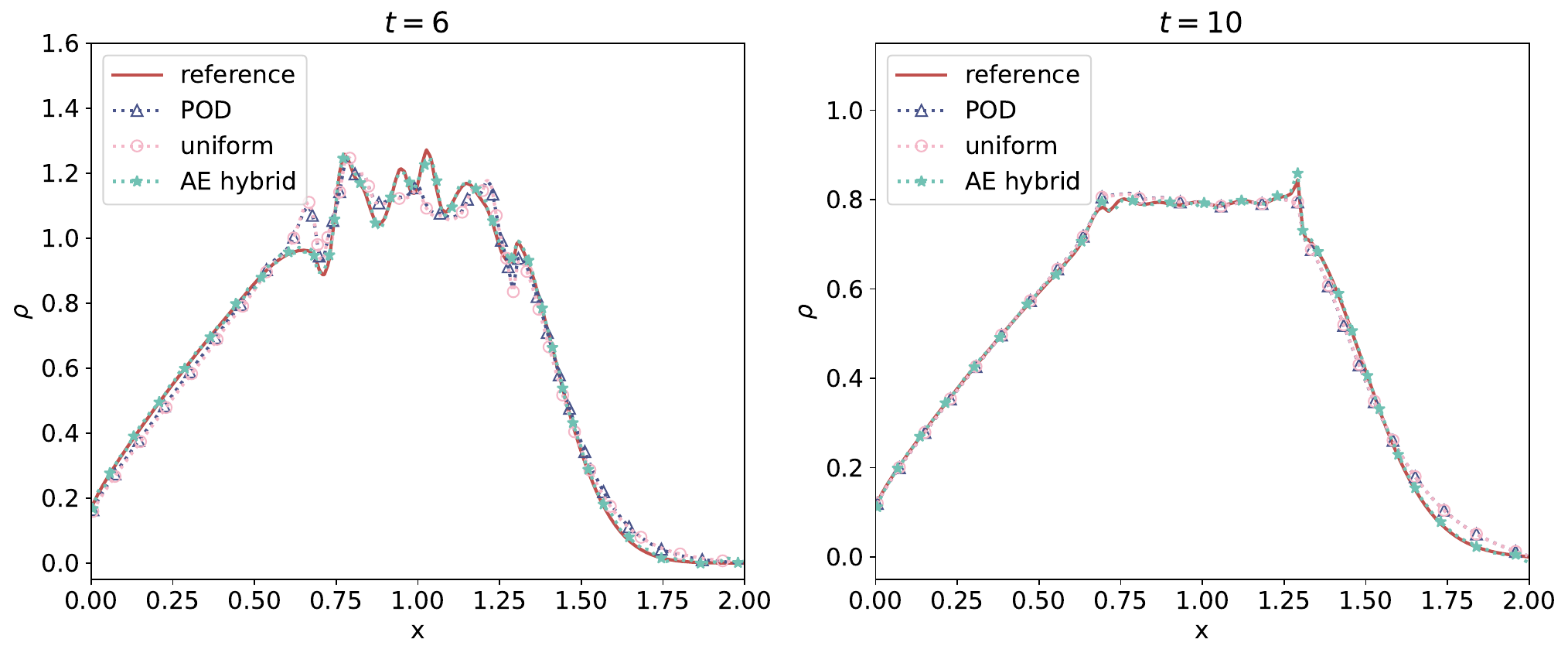}
    		\end{minipage}
    	}%
	\caption{$\rho$ in Example 2 when $\mu=136.10$ and $\varepsilon_{\text{POD}} = 0.1$. (a) $t=3$ and $t=4$; (b) $t=6$ and $t=10$.}
	\label{fig:eg2_hybridrho}
\end{figure}

The details of the hyperparameter settings in this example can be found in \ref{apx:ae}. The total number of trainable parameters in this convolutional autoencoder is $63972$.
Table \ref{table:eg2_hybrid} presents the performance of the autoencoder hybrid method in comparison with classical POD and uniform partitioning POD under a tolerance of $\varepsilon_{\text{POD}}=0.1$.
 By examining the overall online time and relative errors, it is evident that the hybrid approach outperforms both the classical POD and the uniform partitioning POD in terms of accuracy and online efficiency, while using a relatively loose POD tolerance. Specifically, compared to classical POD, uniform-$2$ improves $E_f$ by $12.54$\%, $E_\rho$ by $16.01$\%, and accelerate the online stage by $48.00$\%. The AE hybrid method further enhances these improvements, achieving reductions of $34.96$\% in $E_f$, $23.17$\% in $E_\rho$, and $88.00$\% in online computational time. Furthermore, Figure \ref{fig:eg2_hybridrho} displays the graphs of $\rho$ reconstructed using these three methods at four sample times $t=3,4,6,10$. An interesting trend emerges: at early times ( $t=3,4$), where the reference solution exhibits high oscillations, the hybrid method introduces slight biases. However, at later times ($t=6,10$), where the reference solution becomes smoother, the hybrid method provides a significantly better approximation than both the classical POD and uniform partitioning POD. We suspect that this arises from the neural network's difficulty in accurately capturing highly oscillatory features, whereas smoother solutions are more amenable to autoencoder-based approximation. Both Table \ref{table:eg2_hybrid} and Figure \ref{fig:eg2_hybridrho} highlight the potential of the hybrid approach in achieving high computational efficiency with better reconstruction accuracy, particularly in regimes where the solution exhibits moderate smoothness.
\begin{table}[htbp]
\centering
\footnotesize
\begin{tabular}{|c|c|c|c|c|}
\hline
method       & reduced space dimension & online time (sec) & $E_f$     & $E_\rho$  \\ \hline
classical POD &    [20]        &   0.25       & $9.01\times10^{-2}$ & $5.87\times10^{-2}$ \\ \hline
uniform-2        & [20, 2]       & 0.13    & $7.88\times10^{-2}$ & $4.93\times10^{-2}$ \\ \hline
\textbf{AE hybrid} & \textbf{[4, 2]}   & \textbf{0.03}  & $\mathbf{5.86\times10^{-2}}$ & $\mathbf{4.51\times10^{-2}}$ \\ \hline
\end{tabular}
\caption{Comparison of autoencoder hybrid method with classical POD and uniform time partitioning POD in Example 2, $\varepsilon_{\text{POD}}=0.1$.}
\label{table:eg2_hybrid}
\end{table}

\subsection{Example 3: pin-cell problem}
Finally, we consider a 2-dimensional pin-cell problem on $\textbf{x} = (x_1,x_2)\in\Omega_{\textbf{x}}=[-1,1]^2$ with zero source term and zero inflow boundary conditions. Prior to introduce the rest settings, let us define the inner and outer region of $\Omega_{\textbf{x}}$:
\begin{align*}
    \Omega_{inner} &= \{(x_1,x_1):|x_1|\leq0.5 \text{ and }|x_2|\leq0.5\},\\
    \Omega_{outer} &= \{(x_1,x_1):|x_1|>0.5 \text{ or }|x_2|>0.5\}.
\end{align*}
Then
\begin{align*}
    \sigma_a(x)&=\begin{cases}
0\hspace{17pt} \textbf{x}\in\Omega_{outer},\\
\mu_a \hspace{11pt} \textbf{x}\in\Omega_{inner},
\end{cases}\\ \sigma_s(x)&=\begin{cases}
100,\hspace{6pt} \textbf{x}\in\Omega_{outer},\\
\mu_s, \hspace{11pt} \textbf{x}\in\Omega_{inner},
\end{cases}\\
f(\textbf{x},\textbf{v},0) &= \frac{10}{2\pi}e^{-5000(x_1^2+x_2^2)}.
\end{align*}
In this problem, the inner region has the scattering and absorption coefficients given by $(\mu_s,\mu_a) \in [0.02,0.5]^2$. This means that the scattering effect in the outer region is $200$ to $5000$ times as strong as that in the inner region, resulting in significant multiscale effects.

We use a $40\times40$ uniform mesh on the spatial domain $\Omega_{\textbf{x}}=[-1,1]^2$ and apply a Chebyshev-Legendre discretization on $S^2$ 
This configuration yields $n_x = 6400$ spatial degrees of freedom and $n_v = 300$ velocity degrees of freedom. The time step size is fixed to $\Delta t=0.01$. To generate our dataset, we collect snapshots at $n_t=20$ uniformly spaced time points in the interval $[0.1,2]$, with
$n_p=81$ uniformly sampled parameter values in $[0.02,0.5]^2$. These parameter pairs are arranged in a grid pattern as $\{(\mu_s^{(1)},\mu_a^{(1)}), \dots, (\mu_s^{(1)},\mu_a^{(9)}), \dots, (\mu_s^{(9)},\mu_a^{(9)})\}$. The snapshots for the first $64$ parameters constitute the training set, while the remaining $17$ are reserved for out-of-distribution testing. This differs from previous examples by focusing on generalization to purely unseen parameter configurations. In the online parameter extrapolation, we adopt the RBF interpolation scheme with a quintic kernel.

\begin{figure}[htbp]
	\centering
       \subfigure[]{
		\begin{minipage}{1\textwidth}
			\includegraphics[width=1\textwidth]{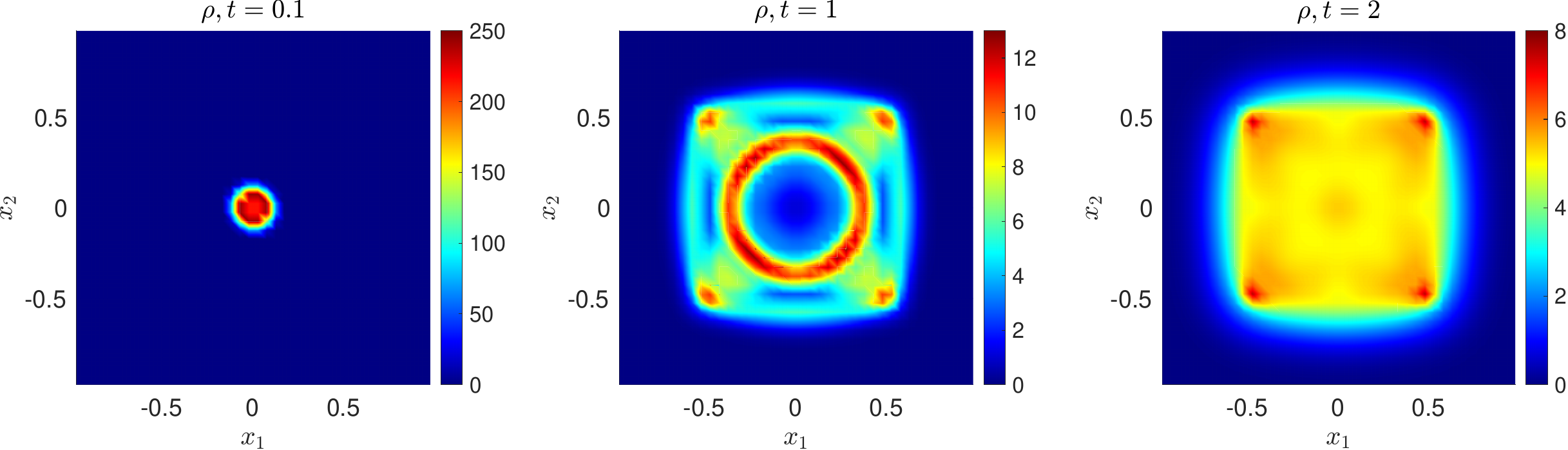}
		\end{minipage}
	}%
    
    	\subfigure[]{
    		\begin{minipage}{1\textwidth}
   		 	\includegraphics[width=1\textwidth]{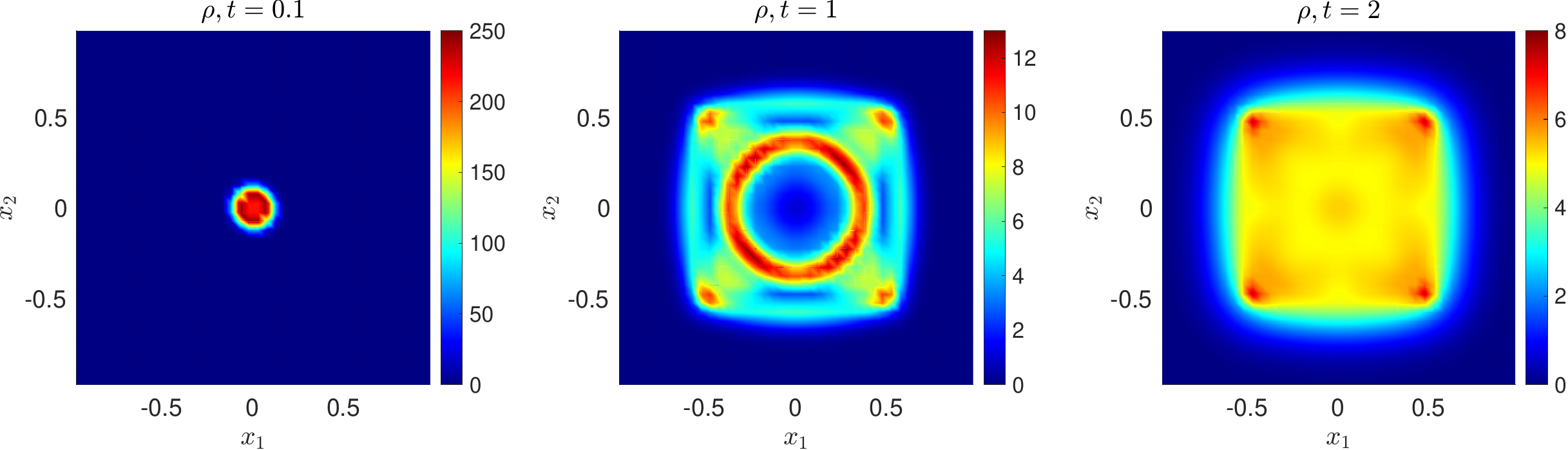}
    		\end{minipage}
    	}%

        \subfigure[]{
    		\begin{minipage}{1\textwidth}
   		 	\includegraphics[width=1\textwidth]{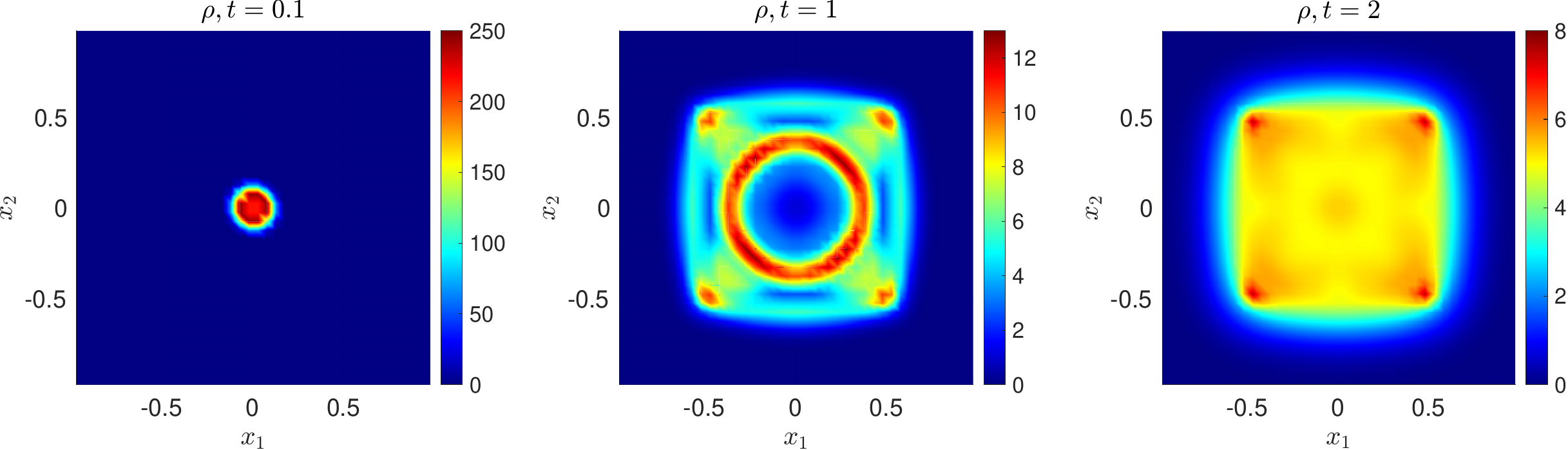}
    		\end{minipage}
    	}%
	\caption{$\rho$ in Example 3: (a) Reference solution; (b) classical POD when $\varepsilon_{\text{POD}} = 10^{-4}$; (c) adaptive method with $\hat{r}_{\max}=20,\hat{r}_{\min}=10$ when $\varepsilon_{\text{POD}} = 10^{-4}$.}
	\label{fig:eg3_rho}
\end{figure}

\textbf{Comparison between uniform and adaptive partitioning.}
We first investigate the performance of the uniform and adaptive time partitioning PODs in 2D case. For comparison, we select classical POD and uniform partitioning into $4$ time intervals as the baseline methods, against which we evaluate the adaptive method with $\hat{r}_M=20$ and $\hat{r}_m=10$. The relative errors, $E_f$, $E_{\rho}$, along with the online CPU time are presented in Table \ref{table:eg3_uni_adap} across $17$ testing parameters. The results demonstrate that time partitioning enhances the accuracy of the solution. An interesting observation is that, as the POD tolerance decreases, the improvement in precision becomes less pronounced, likely due to the dominance of extrapolation error at lower tolerances. Figure \ref{fig:eg3_rho} (c) displays the reconstructed $\rho$ using adaptive method for the test parameter $(\mu_s,\mu_a)=(0.44, 0.08)$ at three sampling times, $t=0.1,1,2$. The predicted solutions match well with the reference solution, demonstrating the effectiveness of the proposed approach.

\begin{table}[htbp]
\centering
\begin{tabular}{|c|c|c|c|}
\hline
         & $E_f$    & $E_\rho$ & online time (sec) \\ \hline
POD      & $1.09\times10^{-2}$ & $3.61\times10^{-3}$ & 0.50             \\ \hline
uniform-4  & $\mathbf{1.08\times10^{-2}}$ & $\mathbf{3.51\times10^{-3}}$ & 0.50             \\ \hline
adaptive & $\mathbf{1.08\times10^{-2}}$ & $3.59\times10^{-3}$ & \textbf{0.40}\\ \hline
\end{tabular}
\caption{Comparison of relative errors and online time of POD, uniform partitioning and adaptive partitioning methods in Example 3}
\label{table:eg3_uni_adap}
\end{table}

\begin{figure}[htbp]
\centering
	\includegraphics[width=0.5\textwidth]{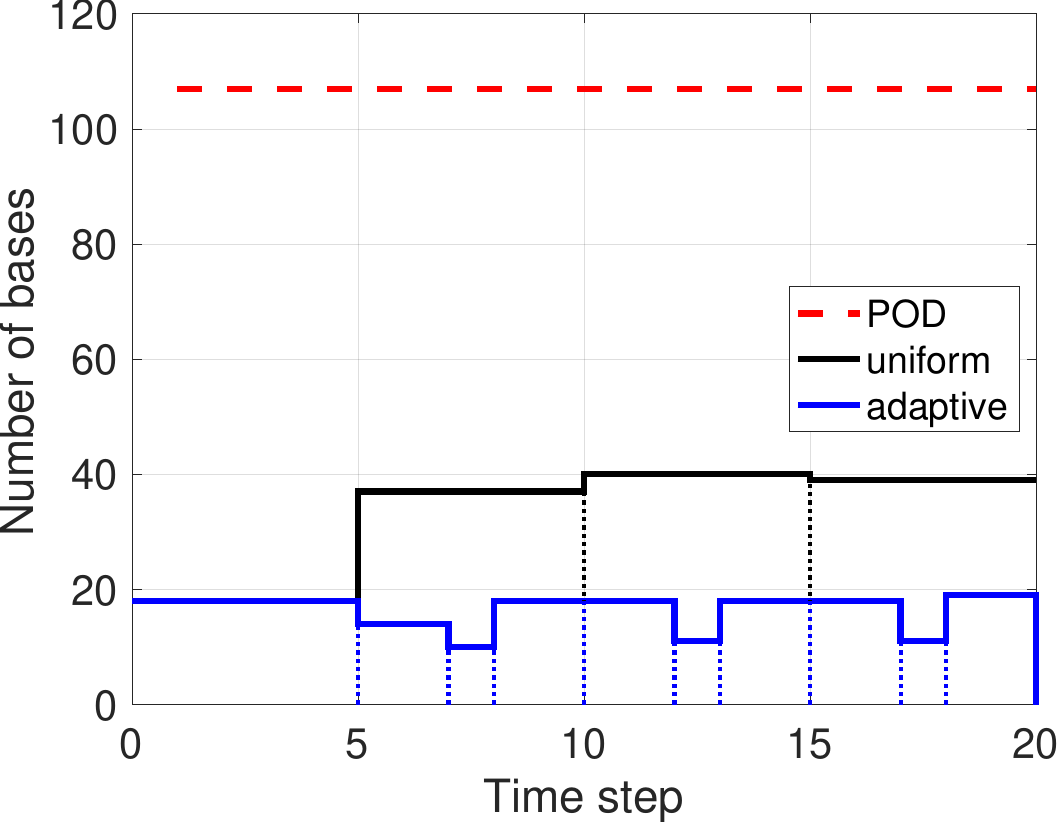}
	\caption{Illustration of number of basis functions on each interval when using different methods in Example 3.}
	\label{fig:eg3_dim}
\end{figure}

 Figure \ref{fig:eg3_dim} illustrates the number of bases required in each time interval using classical POD, uniform-4 and adaptive partitioning POD. We observe that, following uniform time partitioning, the maximum and average number of bases is significantly reduced, which facilitates the computational time complexity in online stage.
Uniform-4 partitioning, however, proves inadequate in the later three time intervals, and the adaptive method controls the number of bases in these intervals below $\hat{r}_M$ through finer time partition.

\textbf{Strength and weakness of the hybridized strategy.}
Since the quantity $\rho_h$ is experimentally observable, rather than the discrete distribution function $f_h$, we focus our analysis on the simulation performance of $\rho_h$ throughout this section. As illustrated in Figure \ref{fig:eg3_dim}, the number of bases in the time interval $[0.5,2]$ is notably higher than that in $[0,0.5]$. Hence, we then consider applying the autoencoder hybrid method here with input parameters $\tau_{\min}=0.5$.

The detailed hyperparameter and training settings can be found in \ref{apx:ae}. The whole convolutional antoencoder under this setting has $1052237$ parameters. Given the training snapshots $\{\rho_h(t_i;\mu_j): i=1,\dots,n_t; j=1,\dots, n_p\}$, now we determine the trainable parameters $\theta$ via minimizing the following loss function:
\begin{align*}
    \mathcal{L} = \frac{1}{n_t n_p}\sum_{i=1}^{n_t}\sum_{j=1}^{n_p} \|\rho_h(t_i, \mu_j)-\mathcal{D}_\theta\circ\mathcal{E}_\theta(\rho_h(t_i, \mu_j))\|^2_{l^2}
\end{align*}

\begin{figure}[htbp]
	\centering
       \subfigure[]{
		\begin{minipage}{1\textwidth}
			\includegraphics[width=1\textwidth]{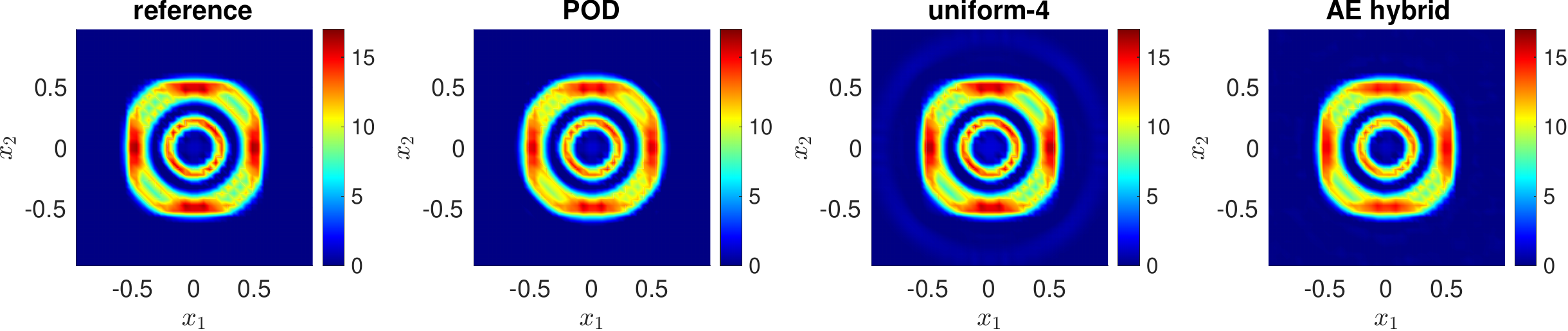}
		\end{minipage}
	}%
    
    	\subfigure[]{
    		\begin{minipage}{1\textwidth}
   		 	\includegraphics[width=1\textwidth]{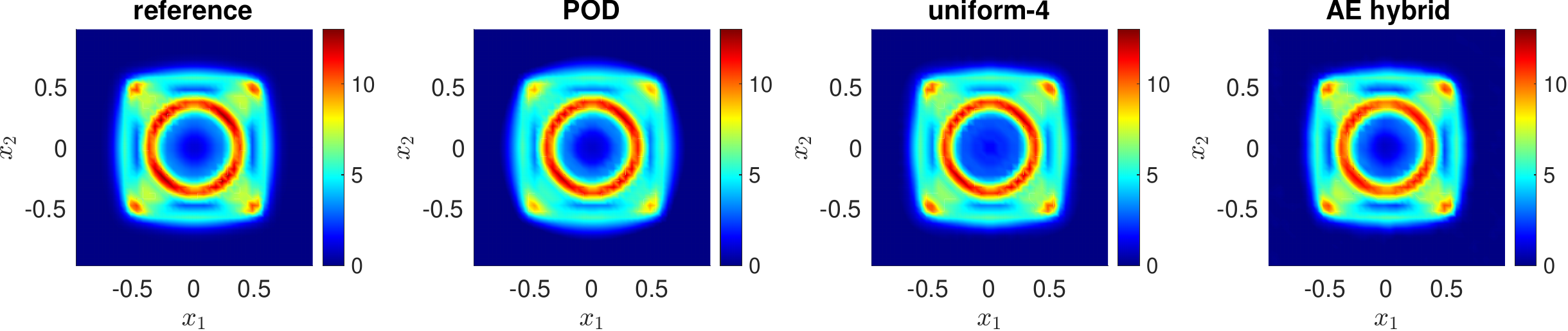}
    		\end{minipage}
    	}%
        
        \subfigure[]{
		\begin{minipage}{1\textwidth}
			\includegraphics[width=1\textwidth]{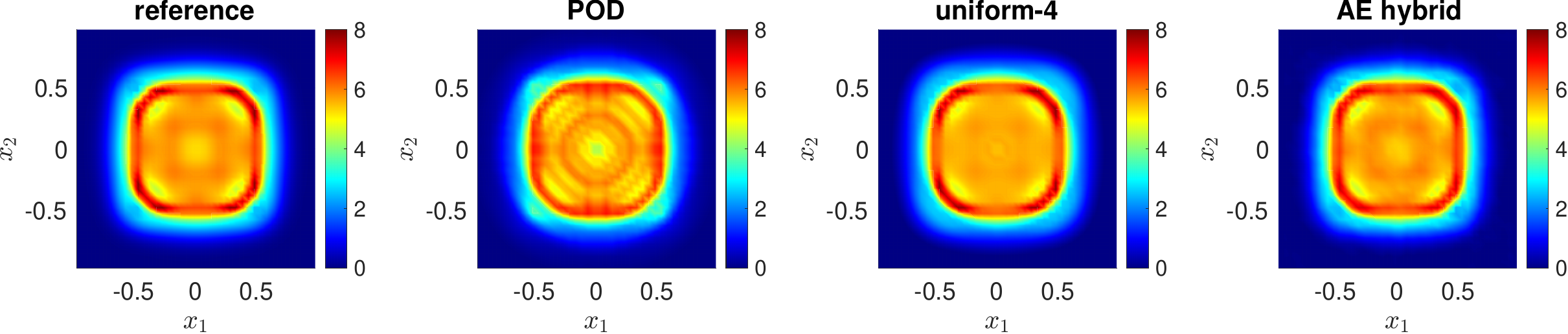}
		\end{minipage}
	}%
    
    	\subfigure[]{
    		\begin{minipage}{1\textwidth}
   		 	\includegraphics[width=1\textwidth]{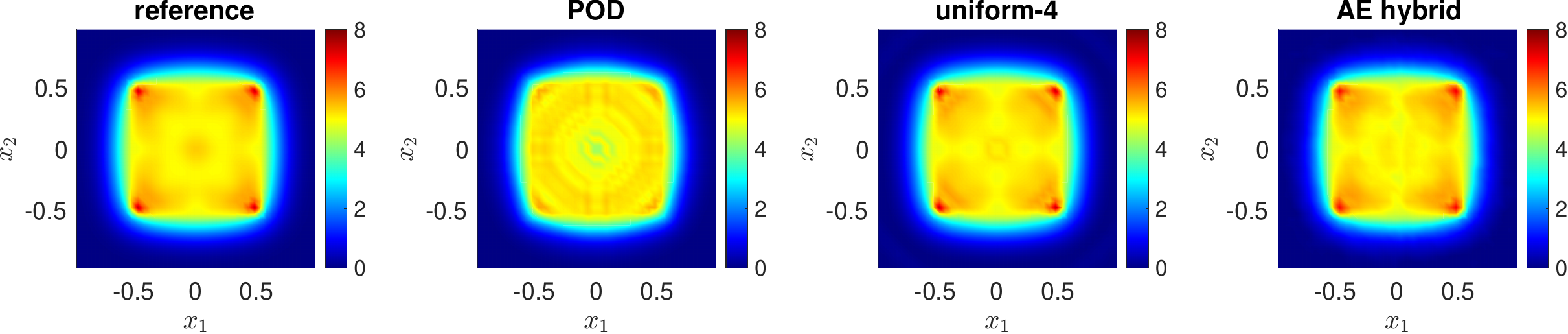}
    		\end{minipage}
    	}%
	\caption{$\rho$ in Example 3 reconstructed by different methods when $\varepsilon_{\text{POD}} = 0.07$. (a) $t=0.6$; (b) $t=1.0$; (c) $t=1.6$; (d) $t=2.0$.}
	\label{fig:eg3_hybridrho}
\end{figure}

Table \ref{table:eg3_hybrid} compares the performance of the hybrid approach against both classical POD and uniform-$4$ at POD tolerance levels $\varepsilon_{\text{POD}}=7\times10^{-2}$ and $10^{-2}$. When $\varepsilon_{\text{POD}}=7\times10^{-2}$, we choose $\hat{r}_M=5$ in the hybridized strategy, while this parameter is chosen as $\hat{r}_M=10$ when $\varepsilon_{\text{POD}}=10^{-2}$.
 The results demonstrate that the hybrid method significantly reduces online computation time for both POD tolerance levels. At $\varepsilon_{\text{POD}}=7\times10^{-2}$, the hybrid T-P $4$ achieves superior accuracy compared to other methods, while maintaining comparable precision to classical and uniform partitioning PODs when $\varepsilon_{\text{POD}}=10^{-2}$. Figure \ref{fig:eg3_hybridrho} illustrates $\rho$ values at four times $t=0.6,1,1.6, 2.0$ with parameter $(\mu_s,\mu_a)=(0.44,0.08)$ across different methods. The figures indicate that the solutions reconstructed by AE hybrid approach consistently show the close agreement with reference values.
\begin{table}[htbp]
  \footnotesize
  \centering
\subfigure[$\varepsilon_{\text{POD}}=0.07$]{
\begin{tabular}{|c|c|c|c|}
\hline
method       & number of basis functions & online time (sec) & $E_\rho$  \\ \hline
classical POD &           [11]        &  $5.00\times10^{-3}$   &$1.20\times10^{-1}$ \\ \hline
uniform-4        & [5, 10, 9, 6]        & $9.12\times10^{-3}$   & $\mathbf{5.16\times10^{-2}}$ \\ \hline
AE hybrid & \textbf{[5, 4, 4, 4]} & $\mathbf{1.66\times10^{-3}}$ & $5.43\times10^{-2}$ \\ \hline
\end{tabular}
  }

\subfigure[$\varepsilon_{\text{POD}}=0.01$]{
\begin{tabular}{|c|c|c|c|}
\hline
method       & reduced space dimension & online time (sec) & $E_\rho$  \\ \hline
classical POD &             [26]        &  $9.99\times10^{-3}$   &$2.41\times10^{-2}$ \\ \hline
uniform-4        & [6, 16, 16, 13]          & $1.66\times10^{-2}$  & $\mathbf{1.50\times10^{-2}}$ \\ \hline
\textbf{AE hybrid} & \textbf{[6, 4, 4, 4]}& $\mathbf{1.85\times10^{-3}}$ & $5.28\times10^{-2}$ \\\hline
\end{tabular}
}
  \caption{Comparison of autoencoder hybrid method with classical POD and uniform time partitioning POD (4 intervals) in Example 3.}
\label{table:eg3_hybrid}
\end{table}

\section{Conclusions}\label{sec:conclusion}
This work introduces computationally efficient nonlinear ROMs based on time partitioning techniques for RTE with multiscale properties. We introduce a novel adaptive time partitioning strategy with both refinement and coarsening abilities, enabling robust use of arbitrary initial partition without concern for superfluous local refinement. An equilibrium detection step is also designed to further improve the online efficiency for long-time simulations.
In addition to the piecewise linear in time ROM, we develop a hybrid framework that adaptively identifies intervals where local linear ROM is insufficient and applies an autoencoder-based nonlinear ROM in these regions. This novel hybrid approach, first time in the literature, combines the advantages of both piecewise linear ROM and neural network-based ROM while simultaneously reducing the computational overhead associated with the training of neural network. 

The efficacy and computational efficiency of our proposed ROMs are further validated through multiple numerical experiments involving kinetic transport equations that encompass both parabolic and hyperbolic scales. Numerical tests demonstrate that the proposed nonlinear ROMs achieve superior performance over traditional linear ROM approaches in terms of accuracy and online computational efficiency. 
Furthermore, the adaptive approach enables precise control over the maximum and minimum number of bases in each time interval, allowing us to tailor the online computational cost to specific application requirements. Notably, the autoencoder hybrid approach successfully identifies and exploits ultra-low-rank structures of solution manifolds, facilitating exceptionally fast online computation.

Several interesting directions for future research center on the theoretical analysis of the proposed ROMs. A primary objective is to establish theoretical foundations of the relation between time partitioning and the decay rate of the Kolmogorov $n$-width \eqref{eq:kol_n_width} in each interval. While \cite{Bachmayr_Cohen_2017} have analyzed the effects of uniform parameter space partitioning on Kolmogorov n-width for elliptic equations, the understanding of the effects of time partitioning for time-dependent equations is still unexplored. Furthermore, the investigation of autoencoder's approximation and generalization errors, particularly in the context of time partitioning, presents another compelling research direction. Despite the previous exploration of approximation and generalization errors of autoencoders in operator learning \cite{AE_generalizationerror}, such theoretical analysis is currently hindered by the lack of analysis of Kolmogorov $n$-width for RTEs, which poses a fundamental challenge to developing a comprehensive theoretical framework for autoencoder-based ROMs in this context. Moreover, more efficient non-intrusive online strategy can be developed.

\section*{Acknowledgments}
T.J. gratefully acknowledges Chutian Huang (HKUST) for providing high performance computational resources funded by her HKPFS. The work of Z.P. was supported by the Early Career Scheme 26302724. The work of Y.X. was supported by the Project of Hetao Shenzhen-HKUST Innovation Cooperation Zone HZQB-KCZYB-2020083. T.J. and Y.X. would like to thank HKUST Fok Ying Tung Research Institute and National Supercomputing Center in Guangzhou Nansha Sub-center for providing high performance computational resources. 
\appendix
\section{Proof of Theorem \ref{thm:error_frob}}\label{sec:proof_error_frob}
\begin{proof} From the definition of the Frobenius norm, we have
\begin{align*}
    \|S-(S_{r_1}|\dots|S_{r_k})\|_F^2 = \sum_{j=1}^k \|S_j - S_{r_j}\|_F^2\leq \varepsilon^2 \sum_{j=1}^k \|S_j\|_F^2 = \varepsilon^2 \|S\|_F^2.
\end{align*}
\end{proof}

\section{Detailed set-up for architectures and training of autoencoders in various examples}\label{apx:ae}
We summarize the detailed settings of the architecture and training of autoencoders for Example 1,2,3 in Table \ref{tab:eg1_ae_architecture}, \ref{tab:eg2_ae_architecture}, \ref{tab:eg3_ae_architecture}, respectively. 

\begin{table}[htbp]
\centering
\begin{tabular}{ll}
\toprule
\textbf{Encoder} & 
\begin{tabular}[t]{@{}l@{}}
• $4$ Conv1d layers with input-output number of channels\\
\quad - $(12,24)$ \\
\quad - $(24,24)$ \\
\quad - $(24,24)$ \\
\quad - $(24,12)$ \\
• Kernel size: $10$, stride: $2$, zero-padding, $n_v$ channel \\
• Final fully connected layer: $n_{x}/16$ neurons \\
• Latent dimension: $r=8$
\end{tabular} \\
\hline
\textbf{Decoder} & Mirror symmetry with encoder \\
\midrule
\textbf{Training} & 
\begin{tabular}[t]{@{}l@{}}
• Optimizer: Adam (initial learning rate, $10^{-3}$) \\
• Learning rate scheduler: plateau detection \\
\quad - Patience: $5$ epochs \\
\quad - Decay rate: $0.25$ \\
• Total epochs: $5000$\\
• Loss: MSE 
\end{tabular} \\
\bottomrule
\end{tabular}
\caption{Autoencoder architecture and training details for Example 1}
\label{tab:eg1_ae_architecture}
\end{table}

\begin{table}[htbp]
\centering
\begin{tabular}{ll}
\toprule
\textbf{Encoder} & 
\begin{tabular}[t]{@{}l@{}}
• $4$ Conv1d layers with input-output number of channels\\
\quad - $(16,32)$ \\
\quad - $(32,32)$ \\
\quad - $(32,32)$ \\
\quad - $(32,16)$ \\
• Kernel size: $10$, stride: $2$, zero-padding, $n_v$ channel \\
• Final fully connected layer: $n_x/16$ neurons \\
• Latent dimension: $r=4$
\end{tabular} \\
\hline
\textbf{Decoder} & Mirror symmetry with encoder \\
\midrule
\textbf{Training} & 
\begin{tabular}[t]{@{}l@{}}
• Optimizer: Adam (initial learning rate: $2\times10^{-4}$; weight decay: $10^{-7}$) \\
• First $1000$ epochs: \\
\quad Learning rate scheduler: \\
\quad\quad - Step size: $100$ epochs \\
\quad\quad - Decay rate: $0.8$ \\
• Later $29000$ epochs: \\
\quad Learning rate: $10^{-5}$\\
• Total epochs: $30000$\\
• Loss: MSE 
\end{tabular} \\
\bottomrule
\end{tabular}
\caption{Autoencoder architecture and training details for Example 2}
\label{tab:eg2_ae_architecture}
\end{table}

\begin{table}[htbp]
\centering
\begin{tabular}{ll}
\toprule
\textbf{Encoder} & 
\begin{tabular}[t]{@{}l@{}}
• $4$ Conv2d layers with input-output number of channels\\
\quad - $(1,4)$ \\
\quad - $(4,16)$ \\
\quad - $(16,64)$ \\
\quad - $(64,64)$ \\
• Kernel size: $10$, stride: $2$, zero-padding, $1$ channel \\
• Final fully connected layer: $n_\mathbf{x}n_\mathbf{y}/4$ neurons \\
• Latent dimension: $r=4$
\end{tabular} \\
\hline
\textbf{Decoder} & Mirror symmetry with encoder \\
\midrule
\textbf{Training} & 
\begin{tabular}[t]{@{}l@{}}
• Optimizer: Adam (learning rate: $10^{-5}$; weight decay: $10^{-7}$) \\
• Total epochs: $40000$\\
• Loss: MSE 
\end{tabular} \\
\bottomrule
\end{tabular}
\caption{Autoencoder architecture and training details for Example 3}
\label{tab:eg3_ae_architecture}
\end{table}



\bibliographystyle{elsarticle-num-names}
\bibliography{references}

\end{document}